\documentclass[11pt]{amsart}
\include{std}
\usepackage[all]{xy}

\usepackage{psboxit}
\usepackage{ifthen}
\newdimen\minCDarrowwidth
\usepackage{amsmath}
\usepackage{amsfonts}
\usepackage{amssymb}
\usepackage{amscd}
\usepackage{verbatim}

\newcommand{\marginlabel}[1]%
 {\mbox{}\marginpar{\raggedleft\hspace{0pt}\bfseries\sf#1}}

\def\AA{{\mathbb A}}
\def\CC{{\mathbb C}}
\def\PP{{\mathbb P}}
\def\ZZ{{\mathbb Z}}

\def\cB{\mathcal{B}}
\def\cC{\mathcal{C}}
\def\cI{\mathcal{I}}
\def\cJ{\mathcal{J}}
\def\cS{\mathcal{S}}
\def\cX{\mathcal{X}}
\def\cN{\mathcal{N}}

\def\cO{\mathcal{O}}
\def\cM{\mathcal{M}}
\def\cP{\mathcal{P}}
\def\cH{\mathcal{H}}
\def\cU{\mathcal{U}}
\def\si{\sigma}
\def\ts{2\sqrt[5]{\frac{4}{3}}\mu}
\def\Im{\mbox{ Im }}

\DeclareMathOperator{\Ker}{Ker}

\DeclareMathOperator{\rank}{rank}
\DeclareMathOperator{\Spec}{Spec}

\newtheorem{lemma}{Lemma}[section]
\newtheorem{theorem}[lemma]{Theorem}
\newtheorem{corollary}[lemma]{Corollary}
\newtheorem{proposition}[lemma]{Proposition}

\theoremstyle{definition}

\newtheorem{remark}[lemma]{Remark}

\newtheorem{notation}{Notation}

\numberwithin{equation}{section}

\newcommand{\bean}{\begin{eqnarray}}
\newcommand{\eean}{\end{eqnarray}}
\newcommand{\bea}{\begin{eqnarray*}}   
\newcommand{\eea}{\end{eqnarray*}}
\newcommand{\be}{\begin{displaymath}}
\newcommand{\ee}{\end{displaymath}}

\begin{document}

\title{Degree 1 Curves in the Dwork Pencil and the Mirror Quintic} 

\author[Anca Musta\c{t}\v{a}]{Anca~Musta\c{t}\v{a}}
\address{Department of Mathematics, University of British Columbia,
Room 121, 1984 Mathematics Road, Vancouver, B.C., Canada, V6T1Z2}
\email{{\tt amustata@math.ubc.edu}}

\date{\today}


\begin{abstract}
We give a description of the relative Hilbert scheme
of lines in the Dwork pencil of quintic threefolds. 
We describe the corresponding relative Hilbert scheme
associated to the mirror family of quintic threefolds.

\end{abstract}

\maketitle
\bigskip

\section*{Introduction}

\subsection*{0.1. The relative Hilbert scheme of lines in the Dwork family of quintics}

 This paper takes a look at a very well known family of Calabi-Yau threefolds: the family $\rho:\cX \to \CC $ whose fiber $X_t$ over $t \in \CC$ is given in $\PP^4$ by the equation $(F_t=0) $, where
\[F_t(z_0:z_1:z_2:z_3:z_4)=z_0^5+z_1^5+z_2^5+z_3^5+z_4^5-5tz_0z_1z_2z_3z_4. \]

 This family is known as the Dwork pencil of quintics.

 We will study the relative Hilbert Scheme of lines in the family $\cX$. Let $\pi: \cH \to \CC $ denote this relative scheme.
We show that the structure of $\cH$ is 
\[ \cH = \cS \cup \bigcup_{i=1}^{375} Z_i \]
where $\cS$ is a smooth irreducible surface, proper over $\CC$ and for each $i$,\\
$\pi_{| Z_i} :Z_i \to \CC$ is an isomorphism. 

 The general fiber of $\pi_{| \cS}: \cS \to \CC$ is a smooth proper curve consisting of two isomorphic components, each of genus 626. 

 The Stein factorization of $\pi_{| \cS}$ is
\bea \begin{CD} \cS @> {\pi_2} >> C_2 @> {\beta} >> \CC \end{CD} \eea
 where $\beta$ extends to a double cover $\bar{\beta} : \bar{C_2} \to \bar{\CC}=\PP^1$ branched at six points: $t=0$ and the fifth roots of $\frac{2^7}{3}$.

 Let $\tilde{G}(2,5)$ denote the blow-up of the Grassmannian $G(2,5)$ of lines in $\PP^4$, along the subvarieties parametrizing lines inside the coordinate hyperplanes of $\PP^4$. The closure $\tilde{\cS}$ of $\cS$ in $\tilde{G}(2,5)\times \PP^1$ is a smooth projective surface admitting a 30-1 finite morphism to the Fermat surface
\bea z_0^5+z_1^5+z_2^5+z_3^5=0\eea
 in $\PP^3$.

 The fiber $\tilde{\cS}_{\infty}$ is the stable limit of the family $\cS$ at $t=\infty$ and consists of two connected components $\tilde{\cC}^{\xi}$ and $\tilde{\cC}^{\xi^2}$, which are isomorphic and reducible:
\bea \tilde{\cC}^{\xi}= \bigcup_{i=0}^4 \tilde{\cC}^{\xi}_i \eea
with $\tilde{\cC}^{\xi}_i$ smooth isomorphic curves of genus 76, intersecting pairwise transversely at 25 points - a total of 500 points.

\subsection*{ The ``mirror'' Hilbert scheme}

 The group $$\hat{G}=\{(\mu_0:\mu_1:\mu_2:\mu_3:\mu_4) /\mu_i^5=1\}/\{(\mu:\mu:\mu:\mu:\mu) / \mu^5=1\}$$ acting on $\PP^4$ also acts equivariantly on $\cS$, $C_2$ and $\CC$. The action on $\CC$ is via $\mu = \prod_{j=0}^4 \mu_j$. Putting parameter $w=t^5$ on $\CC/\hat{G}$ we have the Stein factorization
\bea  \begin{CD} (\cS/\hat{G}) @> {\pi_2/\hat{G}} >> (C_2/\hat{G}) @> {\beta/\hat{G}} >> (\CC/\hat{G}) \end{CD} \eea
  of $\psi$ where $\beta/\hat{G}$ extends to a double cover of $\PP^1$ branched at 0 and $\frac{2^7}{3}$.

 The general fiber of $\pi_2/\hat{G}$ is a smooth hyperelliptic irreducible curve of genus 6. 

 The compactification $\tilde{\cS}/\hat{G}$ of $\cS/\hat{G}$ is a smooth, irreducible, rational surface.

 The stable limit at $ w =\infty$ of the family $\cS/\hat{G}$ is 
$$\tilde{\cS}_{\infty}/\hat{G} = (\tilde{\cC}^{\xi}/\hat{G}) \bigcup (\tilde{\cC}^{\xi^2}/\hat{G})$$ where $\tilde{\cC}^{\xi}/\hat{G}$ is the union of 5 smooth rational curves intersecting pairwise transversely.

\subsection*{0.2. History of problem}

 For positive integers $m$ and $n$, the Fermat Variety $X^n_m$ is the hypersurface in $\PP^n$ given by equation
$$ z_0^m+z_1^m+...+z_n^m=0.$$
 The Fermat varieties have long been subject to the attention of
mathematicians: the Fermat curves have been intensely studied by
number theorists while the Fermat hypersurfaces in general offered
interesting examples in Hodge theory.

 The Dwork pencil is known today mostly for its quotient by the action of the group $\hat{G}$ defined in part (0.1), which after desingularization gives the mirror to the family of quintic threefolds. Mirror symmetry arose in 
1991, when Candelas, de la Ossa, Green and Parkes used that principle to predict invariants very closely related to the expected numbers of rational curves of any degree on a quintic threefold. In the same year, A.Albano and S. Katz gave a first description of the Hilbert scheme of lines on the Fermat quintic threefold $X_0$ (see \cite{ak1}). They proved that all the lines in $X_0$ are actually contained in 50 cones over some plane quintics. They also studied the deformation theory of these lines for general deformations of $X_0$.

Then they initiated the study of the relative Hilbert scheme of the Dwork pencil of quintics (see \cite{ak2}). Their work built on  B. van Geemen's discovery of a set of 5000 special lines on each member $X_t$ of the family - since this number exceeded the known invariant of 2875 lines, this established the Dwork pencil as the only known nontrivial family of quintics having each a continuous family of rational curves.

 Also in \cite{ak2}, A.Albano and S.Katz presented another proof of the fact that the Hilbert scheme $\cH$ of the Dwork pencil is two dimensional. Their method consists of describing the lines contained in some $X_t$ as lines in $\PP^4$ meeting each component of $B$, where $B$ is the base locus of the family $\cX$. Since $B=\bigcup_{i=0}^{4} B_i$ with
$$B_i=\{(z_0:z_1:z_2:z_3:z_4) \in X_0  / z_i=0, \sum z_j^5=0\}$$
the Fermat surfaces in the hyperplanes $H_i= \{z_i=0\} $, they show that the lines meeting each component of $B$ are parametrized by a complete intersection surface $I$ in the Grassmannian $G(2,5)$, given in Pl\"{u}cker coordinates by the equations:
\bea \sum_{j=0}^{4} p_{ij}^5 = 0 \eea
 for $i \in \{1,...,4\}$. However, this surface has a rather large number of irreducible components, only one of which will correspond to the Hilbert scheme $\cH$. For example, for each $B_{ij}= B_i\cap B_j$, there are components of $I$ made of lines meeting $B_{ij}$ and the other components $B_k$ of the base locus. In this paper we isolate the Hilbert scheme from these other components and describe it more precisely.

 This approach of A.Albano and S.Katz  has proven fruitful in studying the closure of $\cH$ in $G(2,5)\times \PP^1$ and its desingularization, as well as in computing the class of these lines in the cohomology ring of $G(2,5)$ - as seen in sections 8 and 9, respectively.

 On the other hand, for $d\geq 2$, little is known today about the degree d rational curves lying in the Dwork pencil and their images in the mirror family. Our hope is that some of the proprieties of the degree 1 Hilbert scheme described in this paper might prove useful in studying this problem.

\subsection*{0.3. Method of proof} 
  
 In analyzing the degree 1 relative Hilbert scheme of $\rho$, a special role is played by those lines in $X_t$, $t\not=0$, which contain a point two of whose coordinates are 0. We call these lines van Geemen lines, due to their use by Bert van Geemen in \cite{ak1} to give the first proof of the fact that the relative dimension of $\cH=Hilb_1(\cX/\CC) $ is 1. It turns out that the local behavior of $\cH$ at the Van Geemen lines is geometrically understandable. On the other hand, if we remove the lines in $X_0$ and the van Geemen lines from the universal line $\cU$, the resulting quasi-projective scheme $\cU'$ admits a simple and uniquely symmetric algebraic presentation giving a birational map to $\PP^3$.

 The main tool in linking these two viewpoints came from a recent result of Clemens in \cite{c1} for K-trivial threefolds. There, the local analytical Hilbert scheme of a smooth curve is described as a gradient variety.

\subsection*{0.4. Notations}
  Throughout this paper, $\PP^n$ will denote the n-dimensional complex
  projective space, $G(k+1,n+1)$ will denote the Grassmann variety of
  $k$-dimensional projective subspaces of $\PP^n$.

 The Dwork pencil of quintics is denoted by $\cX$ and its relative
 Hilbert scheme of lines by $\cH$. 

 $\hat{G}$ is the group
  $$\{(\mu_0:\mu_1:\mu_2:\mu_3:\mu_4) /\mu_i^5=1\}/\{(\mu:\mu:\mu:\mu:\mu) / \mu^5=1\}$$ 
 and $G$ is its subgroup
 $$\{(\mu_0:\mu_1:\mu_2:\mu_3:\mu_4)
/\mu_i^5=1, \prod_{i=0}^4\mu_i =1\}/\{(\mu:\mu:\mu:\mu:\mu) / \mu^5=1\}$$

\subsection*{0.5. Acknowledgments}
 This paper would not have existed without the generous support of my advisor, Herbert Clemens. I am grateful for his sharing of mathematical ideas, his patience and for the sense of enjoyment that resulted from working with him. I would also like to thank his family for their kindness and hospitality.
 
 I would like to thank Alberto Albano, Sheldon Katz and Bert van
 Geemen, both for setting up this problem, and for their sharing
 insights at different stages of the project. Thank you Gianluca
 Pacienza for useful conversation and encouragement. Mike Stillman, by
 teaching me Macaulay2, has added an exciting aspect to this work. I
 am grateful to Aaron Bertram for his support during these years at the University of Utah.

\newpage

\subsection*{0.6. Contents}

\smallskip
\subsection*{Section 1} Hilbert scheme $\cH_0$ of lines in the Fermat
quintic threefold.
\subsection*{Section 2} The relative Hilbert scheme $\cH$ in a neighborhood of $\cH_0$: the Hilbert scheme as a gradient scheme.
\subsection*{Section 3} The van Geemen lines.
\subsection*{Section 4} 
Algebraic description of the universal line $\cU$ away from $\cU_0$ and the van Geemen lines.
\subsection*{Section 5} 
The surface component $\cS$ of $\cH$.
\subsection*{Section 6} 
Properties of the fiber $( \cS / \hat{G})_{w}$.
\subsection*{Section 7} 
Properties of the fiber $\cS_t$.
\subsection*{Section 8} 
The stable limits $\tilde{\cS}_{\infty }$ and $\tilde{\cS}_{\infty
}/\hat{G}$.
\subsection*{Section 9} 
The class $[\bar{\cS}]$ on $G(2,5)$.

\newpage

\section{ Hilbert scheme $\cH_0$ of lines in the Fermat quintic threefold}

 In \cite{ak1}, $ \cH_{0 red}$ is proven to be the union of 50 Fermat quintic plane curves meeting transversely in pairs at 375 points.  We will start by recalling the description of the Hilbert scheme  $\cH_{0}$ as in \cite{ak1} and \cite{ck}:

\begin{notation} Let $ \mu $ be a complex number such that $\mu ^5 =1 $.\\
 Let $$B_i=\{(z_0:z_1:z_2:z_3:z_4) \in X_0  / z_i=0, \sum z_j^5=0\},$$  for $ i\in \{0,1,2,3,4\}.$

 Let $B_{ij}=B_i\bigcap B_j$ and $B_{ijk}=B_i\bigcap B_j \bigcap B_k$, for $i,j,k \in \{0,1,2,3,4\}.$
\end{notation}

 $\cB = \bigcup_{i=0}^4 B_i$ is the base locus of the pencil of quintic hypersurfaces $X_t$. Each  $B_{ij}$ is a quintic curve in $\PP^2$ and each $B_{ijk}$ is a set of 5 points in $\PP^1$.

 Here $i,j,k,l,h$ will stand for pairwise different indices from 0 to 4.\\
 For each pair $(i,j)$ with $i,j \in \{0,1,2,3,4\} $, the hyperplane $( z_i+\mu z_j=0 )$ intersects $X_0$ in a cone $C_{i,j, \mu}$ over the curve $B_{ij}$, whose vertex is one of the five points of $B_{khl}$. There are 50 such cones. Two such cones  $C_{i,j, \mu}$ and  $C_{k,l, \nu}$ intersect in a unique line iff $i, j, k, l$ are pairwise different. Thus by counting we get 375 lines of intersection: each cone contains 15 such special lines, each special line is at the intersection of exactly 2 cones. The fact that the only lines in the Fermat quintic $X_0$ are those lying in one of these cones is proven in \cite{ak1}. 

 From a local computation it can be seen that the Hilbert scheme $ \cH_{0}$ is not reduced at any of its points ( see \cite{ck}).  We will denote by $l_1$ one of the 375 special points of  $ \cH_{0}$ and by $l_2$ any point  of  $ \cH_{0}$ away from the crossings. Then $ \cH_{0}$ can be described locally analytically as \[ \Spec \CC [[x,y]]/(y^2) \] in a neighborhood of $l_2$ and as  \[ \Spec \CC [[x,y]]/(x^3y^2, x^2y^3) \] in a neighborhood of $l_1$. H.Clemens and H.Kley in \cite{ck} recovered from this local structure the Gromow-Witten invariant for lines in quintic threefolds: namely, the local primary decomposition of the ideal $(x^3y^2, x^2y^3)$: $$ (x^3y^2, x^2y^3) = (x^2)\cap (y^2) \cap (x,y)^5 $$ gives the following computation for the Gromow-Witten invariant:
$$ 50 \cdot 2 \cdot (2g_0-2)+ 375 \cdot 5 = 2875, $$ 
 where $g_0 = 6$ stands for the genus of the plane quintic.

 In the remainder of this section we will answer the question concerning smoothness of the relative Hilbert scheme at the points $l_1$ and $l_2$. We will see that the $l_2$ are smooth points of $\cH$ while the $l_1$ are not. For this purpose we compute the normal bundles $\cN_{l_1\backslash X_0}$, $\cN_{l_1\backslash \cX}$,  $\cN_{l_2\backslash X_0}$, and $\cN_{l_2\backslash \cX}$. Notice also that the lines $\{l_1\}$ are contained in the base locus of the family $(X_t)_{t\in\CC }$, thus the normal bundle  $\cN_{l_1\backslash X_t}$ will also give us information on the structure of $\cH$ around the corresponding points in the fiber $\cH_t$.

 For simplicity, we will consider $l_1$ to be given by $ \phi_{1}: \PP^1 \to \PP^4 $ \[\phi_{1}((\alpha :\beta )):= (\alpha :\beta :-\alpha :-\beta :0) \]
 $l_2$ to be given by the morphism  $ \phi_{2}: \PP^1 \to \PP^4 $ \[\phi_{2}((\alpha :\beta )):= (\alpha :\beta :-\alpha :a\beta :b\beta ) \] with $ab\not= 0$, $ 1+a^5+b^5=0$.

 Using the action of the group $G\times S_5$ we can see that our computations hold for any  $l_1$ and  $l_2$ in $\cH_{0}$. 

 We will use the following lemma, which is a version of lemma 1.8 in \cite{ck} in our particular case:
\begin{lemma}
Let $\cU$ represent the universal family of lines over $\cH$.
\[ \begin{CD}
\cU @> q >>\cX \\
@VV{p}V  \\
\cH
\end{CD} \]
 Let $\cI_{\cH}$ denote the ideal sheaf of $\cH$ in $G(2,5)\times\CC$
 and let $d$ on $G(2,5)\times\CC$ represent the differential in the
 direction of $G(2,5)$. Let $\cI_{\cX}$ denote the ideal sheaf of
 $\cX$ in $\PP^4\times\CC$, $\cI_{\cU}$ the ideal sheaf of $\cU$ in
 $\cH\times_{\CC}\cX$ and  $\cJ_{\cU}$ the ideal sheaf of $\cU$ in
 $\cH\times_{\CC}(\PP^4\times\CC)$. We will denote by $\omega$ the
 relative sheaf of differentials of $\cU$ over $\cH$. 
Then there is a commutative diagram:
 \be \begin{CD}
 R^1p_*(q^*(\cI_{\cX}/\cI_{\cX}^2) \otimes \omega ) @>>>R^1p_*(\cJ_{\cU }/\cJ_{\cU}^2 \otimes \omega) @>>> R^1p_*(\cI_{\cU}/\cI_{\cU}^2\otimes \omega )\\
  @VV{\cong}V  @VV{\cong}V  @VV{\cong}V  \\
\cI_{\cH} /\cI^2_{\cH} @>d>> \Omega^1_{G(2,5)\times\CC/\CC}\otimes \cO_{\cH} @>{j^*}>> \Omega^1_{\cH/\CC}
\end{CD} \ee
 with $j^*$ surjective.

 On fibers over $\{l\} \in \cH_{t}$, the upper morphism in the previous sequence is:
\bea 
H^1(\cN_{X_t\backslash\PP^4}^{\vee}|_l\otimes\cO_l(-2))\to H^1(\cN_{l \backslash \PP^4}^{\vee}\otimes\cO_l(-2))\to H^1(\cN_{l\backslash X_t}^{\vee}\otimes\cO_l(-2))  \eea
which can be thought of, via Serre duality, as 
\bea \begin{CD}
 H^0(\cN_{X_t \backslash \PP^4})|_l)^{\vee} @>>> H^0(\cN_{l \backslash \PP^4})^{\vee} @>>> H^0(\cN_{l \backslash X_t})^{\vee} @>>> 0
\end{CD} \eea
\end{lemma} 
 We may think of $ \Omega^1_{\cH/\CC, \{l\}}$ as the cokernel of the map $\psi_t^{\vee}$, where $\psi_t$ is 
 \bea \begin{CD}
0 @>>> H^0(\cN_{l \backslash X_t}) @>>> H^0(\cN_{l \backslash \PP^4}) @>>> H^0(\cN_{X_t \backslash \PP^4})|_l)\\
        &&  @VV{=}V  @VV{=}V  @VV{=}V\\
0 @>>> H^0(\cN_{l \backslash X_t}) @>>> \bigoplus_{i=1}^3 H^0(\cO_l(1))  @>{\psi_t}>> H^0(\cO_l(5)) 
\end{CD} \eea
 
  Consider  $ U_{01}=G(2,5)\backslash\{p_{01}=0\}$ and  $\{l\} \in U_{01}$ generic given by \\ $ \phi: \PP^1 \to \PP^4 $ \[\phi((\alpha :\beta )):= (\alpha :\beta :x_2\alpha +y_2\beta :x_3\alpha +y_3\beta :x_4\alpha +y_4\beta ) \] the line passing through the points $x=(1:0:x_2:x_3:x_4)$ and \\
$y=(0:1:y_2:y_3:y_4)$, where $(x_2, x_3, x_4, y_2, y_3, y_4)$  are coordinates for $U_{01}$.
 Then
\bea \begin{CD}  \bigoplus_{i=2}^4 H^0(\cO_l(1))  @>{\psi_t}>> H^0(\cO_l(5)) \end{CD} \eea  is \[\psi_t=\left(\phi^*(\frac{\partial F_{t}}{\partial z_2}),\phi^*(\frac{\partial F_{t}}{\partial z_3}),\phi^*(\frac{\partial F_{t}}{\partial z_4})\right). \]
 In particular,  if we consider the canonical bases $(\alpha ,\beta )$
for $H^0(\cO_l(1))$ and $(\alpha^5, \alpha^4\beta,..., \alpha\beta^4,
\beta^5)$ for $H^0(\cO_l(5))$, then 
\bea \begin{CD}  \bigoplus_{i=2}^4 H^0(\cO_{l_2}(1))  @>{\psi_{2t}}>> H^0(\cO_{l_2}(5)), \end{CD} \eea
corresponding to the line $l_2$ embedded by $\phi_2$, is
\bea \psi_{2t}((a_2,b_2),(a_3,b_3),(a_4,b_4))=\sum_{i=2}^4(a_i\alpha +b_i\beta )\left(\phi_2^*(\frac{\partial F_{t}}{\partial z_i})\right)=\eea
\bea =5(a_2\alpha +b_2\beta )(\alpha^4-tab\alpha\beta^3) + \eea
\bea + 5(a_3\alpha +b_3\beta )(a^4\alpha^4+tb\alpha^2\beta^2)+5(a_4\alpha+b_4\beta )(b^4\beta^4+ta\alpha^2\beta^2); \eea
or, in matrix form with respect to the standard bases of $H^0(\cO_{l_2}(1))$ and $H^0(\cO_{l_2}(5))$:
     \bea\psi_{2t}  = 5 \left( \begin{array}{cccccc}
      1 & 0 & 0&-tab&0&0 \\
      0 & 1 & 0&0&-tab&0 \\
      0&0&tb&0&a^4&0\\
      0&0&0&tb&0&a^4\\
      0&0&ta&0&b^4&0\\
      0&0&0&ta&0&b^4 \end{array} \right)\eea  
The determinant of this matrix is  $t^2(a^5+b^5)^2 = t^2$.

 $\psi_{1t}$ may be obtained from the above by setting $a=-1, b=0$. Thus:
\bean h^0(\cN_{l_1\backslash X_0})=\dim\Ker \psi_{10}=2, \eean
\bean h^0(\cN_{l_1\backslash X_t})=\dim\Ker \psi_{1t}=0, \eean for $t \not= 0$. Furthermore
\bean h^0(\cN_{l_2\backslash X_0})=\dim\Ker \psi_{20}=2. \eean

 Similarly, $ h^0(\cN_{l_1\backslash \cX})$, $h^0(\cN_{l_2\backslash \cX})$, can be obtained from the exact sequence
 \be \begin{CD}
H^0(\cN_{l \backslash \cX}) @>>> H^0(\cN_{l \backslash (\PP^4\times \CC )}) @>>> H^0(\cN_{\cX \backslash (\PP^4\times\CC })|_l)\\
       @VV{=}V  @VV{=}V  @VV{=}V\\
H^0(\cN_{l \backslash \cX}) @>>> \bigoplus_{i=1}^3 H^0(\cO_l(1))\oplus H^0(\cO_l) @>{\tilde{\psi_t}}>> H^0(\cO_l(5)) 
\end{CD} \ee
where $\tilde{\psi_t}$ is obtained from $\psi_t$ by adding the row:
\[ ( \begin{array}{cccccc}
      0&0&0&0&0&0 \end{array}), \] 
given by $\phi^*_1(\frac{\partial F_t}{\partial t})$  for $l_1$, and
\[ ( \begin{array}{cccccc}
      0&0&0&-ab&0&0\end{array}) \] given by $\phi^*_2(\frac{\partial
F_t}{\partial t})$ for $l_2$. Thus:
\bean h^0(\cN_{l_1\backslash \cX})=\dim\Ker \tilde{\psi_{1t}}=3 \eean
\bean h^0(\cN_{l_2\backslash \cX})=\dim\Ker \tilde{\psi_{2t}}=2. \eean
 Note that A.Albano and S.Katz have shown in \cite{ak1} that the local dimension of $\cH$ is 2 at each point of $\cH_0$.

We can summarize the  above results in the following:
 \begin{lemma}
With the above notations, we have:
 \[  \begin{array}{ll}
\cN_{l_1\backslash X_0}=\cO_{\PP^1}(1)\oplus\cO_{\PP^1}(-3), & \cN_{l_1\backslash \cX}=\cO_{\PP^1}(1)\oplus\cO_{\PP^1}(-3)\oplus\cO_{\PP^1}\\
\cN_{l_2\backslash X_0}=\cO_{\PP^1}(1)\oplus\cO_{\PP^1}(-3), & \cN_{l_2\backslash \cX}=\cO_{\PP^1}\oplus\cO_{\PP^1}\oplus\cO_{\PP^1}(-2)
\end{array} \] and $\cN_{l_1\backslash X_t}=\cO_{\PP^1}(-1)\oplus\cO_{\PP^1}(-1)$ for $t\not= 0$.

 Thus the relative Hilbert scheme $ \cH$ of the Dwork pencil is smooth at the points of  $ \cH_0$ other than its crossing points. It contains as irreducible components 375 smooth curves corresponding to the base locus of the family $ \cH_{t}$. These curves intersect the rest of $ \cH$ only in the 375 crossing points of $ \cH_{0}$.
\end{lemma}
\begin{proof}
 By Grothendieck's lemma,  $$\cN_{l\backslash X_t}=\cO_{\PP^1}(a)\oplus\cO_{\PP^1}(b)$$  for any of the above mentioned lines $l\subset X_t$, where $a+b=-2$ by adjunction. This, together with formulas (1.2)-(1.6) are enough to establish the form of $\cN_{l\backslash X_t}$. For $\cN_{l\backslash \cX}$, use that \[\cN_{l\backslash \cX}=\cO_{\PP^1}(c)\oplus\cO_{\PP^1}(d)\oplus\cO_{\PP^1}(e)\] with $a\leq c$, $b\leq d$ and $c+d+e=-2$, since  $\cN_{l\backslash \cX}$ is the middle term in the exact sequence
 \be \begin{CD}
0 @>>> \cN_{l \backslash X_t} @>>>\cN_{l \backslash \cX} @>>>\cN_{X_t \backslash \cX}|_l \\
     &&  @VV{=}V  @VV{=}V @VV{=}V\\
0 @>>>  \cO_{\PP^1}(a)\oplus\cO_{\PP^1}(b)  @>>> \cO_{\PP^1}(c)\oplus\cO_{\PP^1}(d)\oplus\cO_{\PP^1}(e) @>>>\cO_{\PP^1} 
\end{CD} \ee
 The local proprieties of $\cH$ are deduced via the identifications of the tangent spaces: $T_{\cH_{t},\{l\} }=H^0(\cN_{l \backslash X_t})$, $T_{\cH,\{l\} }=H^0(\cN_{l \backslash \cX})$.
\end{proof}
 
In fact, some of the above results have already been obtained via different methods in \cite{ak2}.

\bigskip 

\section{The relative Hilbert scheme $\cH$ in a neighborhood of $\cH_0$: the Hilbert scheme as a gradient scheme}

 In \cite{c2}, the relative Hilbert scheme of curves on a family of Calabi-Yau threefolds is constructed locally analytically as a (relative) gradient scheme, more precisely as the zero scheme of the exterior derivative of an analytic function on an analytic family of differentiable curves in the threefolds. In this section we will see how the description of $\cH$ locally as a gradient scheme can give us information on its behavior in the neighborhood of  $\cH_0$. 

 The following lemma will prove useful in this context:

\begin{lemma}          
Let $\cC$ denote one of the 375 curves identified above as components of the relative Hilbert scheme. Let $\cM \subset \cO_{\cC}$ denote the maximal ideal at $t=0$. Then in a neighborhood of 0, the maximal ideal $\cM$ annihilates the stalk of $ \left(\Omega^1_{\cH/\PP^1}\otimes\cO_{\cC}\right)$.
\end{lemma}
\begin{proof}
With the same notations as in the preceding section,
\[H^0\left(U_{01}\cap \cC ,\Omega^1_{\cH/\CC}\otimes\cO_{U_{01} \cap \cC}\right) =\frac{\CC <dx_2,dy_2,dx_3,dy_3,dx_4, dy_4>}{\Im (\psi_{2t})^*}\]
 On $\cC$, $\psi_{2t}$ is given (with respect to the canonical bases) by     \bea\psi_{2t}  = 5 \left( \begin{array}{cccccc}
      1 & 0 & 0&0&0&0 \\
      0 & 1 & 0&0&0&0 \\
      0&0&0&0&1&0\\
      0&0&0&0&0&1\\
      0&0&-t&0&0&0\\
      0&0&0&-t&0&0 \end{array} \right) .\eea  
Then by Lemma 1.1
\[H^0\left(U_{01}\cap \cC ,\Omega^1_{\cH/\CC}\otimes\cO_{U_{01} \cap \cC}\right) \cong \frac{\CC[t]<dx_4, dy_4>}{<tdx_4,tdy_4>} \]
 and therefore $\cM$  annihilates the sections of $\left(\Omega^1_{\cH/\CC}\otimes\cO_{U_{01} \cap \cC}\right)$.
\end{proof}
  The preceding lemma and the fact that locally analytically the Hilbert scheme may be written as a gradient scheme can be used to establish the following proposition:
\begin{proposition}
 Let $l_1$ be one of the crossing points of $\cH_0$. Then there exists an analytical neighborhood $U$ of $l_1$ such that for $t$ close to 0, $\cH_t \cap U $ consists of two smooth disjoint curves and an isolated zero, (which is $\cC \cap \cH_t$). Moreover, $(\cH \backslash \cC )\cap U $ is a smooth irreducible surface. \end{proposition}
\begin{proof}
 As mentioned before, in an analytical neighborhood of $l_1$, the Hilbert scheme $\cH_0$ can be described  as  $ \Spec \CC [[x,y]]/(x^3y^2, x^2y^3) $  (where the coordinates $x,y$ can be thought of as coordinates in an analytic disc of dimension $2= h^0(\cN_{l_1\backslash X_0})$).

   Following \cite{c2}, an easy exercise shows that a potential function $\Phi^0$, whose gradient gives $\cH_0$, can be written, after an appropriate choice of coordinates, as $x^3y^3$. Furthermore, one can extend the function $x^3y^3$ to some analytic function $\Phi(x,y,t)$ such that
\[\Phi(x,y,t)= x^3y^3+ tg(x,y)+t^2h(x,y,t)\] and, in some analytical neighborhood $U$ of $l_1$, $\cH$ is given by the zeroes of the differential of $\Phi$ in the direction of the fiber: $$(\Phi_x=\Phi_y=0).$$ Moreover, we can choose coordinates $x, y$, such that $\cC$ is given by $(x=y=0)$ on U. This means that $g_x(0,0)=g_y(0,0)=0$, i.e. $$g(x,y)= ax^2+2bxy+cy^2+o(3).$$ 
 In this setting,  $\Omega^1_{\cH/\CC}\otimes\cO_{\cC}$ over U has the form
\bea \frac{\CC[t]<dx, dy>}{<\begin{array}{c} (tg_{xx}(0,0)+t^2h_{xx}(0,0,t))dx+(tg_{xy}(0,0)+t^2h_{xy}(0,0,t))dy,\\
(tg_{xy}(0,0)+t^2h_{xy}(0,0,t))dx+(tg_{yy}(0,0)+t^2h_{yy}(0,0,t))dy \end{array} >} \eea
and under these conditions Lemma 2.1. implies that 
\[  \left| \begin{array}{cc}
      g_{xx}(0,0) & g_{xy}(0,0) \\
     g_{xy}(0,0) & g_{yy}(0,0) \end{array} \right|=ac-b^2 \not= 0,\] 
otherwise t would be a nilpotent element of $\Omega^1_{\cH/\PP^1}\otimes\cO_{\cC}$ over U.

  Then \[\Phi_x=3x^2y^3+2t(ax+by+o(2))+t^2h_x(x,y,t),\]
\[\Phi_y=3x^3y^2+2t(bx+cy+o(2))+t^2h_y(x,y,t).\]
In other words,
 \[\Phi_x=3x^2y^3+t(Ax+By),\]
\[\Phi_y=3x^3y^2+t(Cx+Dy),\]
where $A,B,C,D$ are analytical functions in $x,y,t$, with $AD-BC $ invertible in an analytic neighborhood of  $(0,0,0)$.

Notice \[C\Phi_x - A\Phi_y =3x^2y^2(Cy-Ax)+t(CB-AD)y= yk_1(x,y,t), \]
where $k_1(x,y,t) := 3x^2y(Cy-Ax)+t(CB-AD).$
   \[D\Phi_x - B\Phi_y =3x^2y^2(Dy-Bx)-t(CB-AD)x = xk_2(x,y,t), \]
where $k_2(x,y,t) := 3xy^2(Dy-Bx)-t(CB-AD).$ 

 Since $AD-BC\not$ is invertible in an analytic neighborhood of  $(0,0,0)$, it follows that  the Hilbert scheme $\cH_t$ is given in that neighborhood by the ideal \[(yk_1(x,y,t) ,  xk_2(x,y,t) ).  \] 
 On the other hand,  $\cH$ has at least one local component of dimension 2, that cannot be given by $x=0$ or $y=0$ since $x=y=0$ is isolated in $\cH_t$ for small non-zero $t$. So $$k_1=Ek_2$$ for a unit $E$, and the zeroes of $k_1$ give a smooth surface with tangent plane $(t=0)$ at $(0,0,0)$. By inspection on $k_1$ and $k_2$ we get that $E$ must equal 1, $ A=D=0$ and $ B=C$.

 Thus for $U$ small enough, $\cH \cap U$ has irreducible components given by $(x=y=0)$ and $(k_1=0)$, and $\cH_t \cap U $ consists the two smooth disjoint curves $$ xy = \pm t^{1/2} $$
together with the isolated point $x=y=0$.
 \end{proof}
\begin{proposition}
For $t\not= 0$ in an analytical neighborhood of 0, the relative Hilbert scheme $\cH_t$ is a smooth curve, together with 375 isolated points.
\end{proposition}
\begin{proof}
It remains to prove that as it departs from a point ${l_2}$ in the $t$-direction, the non-reduced curve $\cH_0$ splits locally analytically into two non-intersecting analytical open subsets of $\cH_t$, both of which are smooth, 1-dimensional.

 In an analytical neighborhood of $l_2$, the Hilbert scheme $\cH_0$
can be described  as  $ \Spec \CC [[x,y]]/(y^2) $. Here again, as in the previous proposition, the coordinates $x,y$ can be thought of as coordinates in an analytic disc of dimension $2= h^0(\cN_{l_2\backslash X_0})$, which is in accord with the data in  \cite{c2}. Thus one can again extend the function $x^3$ to \[
\Phi(x,y,t)= x^3+ tg(x,y)+t^2h(x,y,t)\] analytical such that, in an analytical neighborhood $U$ of $l_2$, $\cH$ is given by $(\Phi_x=\Phi_y=0)$, with 
 $$ \Phi_x(x,y,t) = 3x^2 +tg_x(x,y)+t^2h_x(x,y), $$
 $$ \Phi_y(x,y,t) = tg_y(x,y)+t^2h_y(x,y). $$
For any $\{l\} \in U\cap \cH$, $\dim T_{\cH, l} = 3 - r(x,y,t)$ where $r(x,y,t)$ is the rank of the following matrix
\bea\left(\begin{array}{ll}
      6x+tg_{xx}+t^2h_{xx}&tg_{xy}+t^2h_{xy}\\
     tg_{xy}+t^2h_{xy}&tg_{yy}+t^2h_{yy} \\
g_x+2th_x+t^2h_{tx}&g_y+2th_y+t^2h_{ty}
      \end{array}\right).\eea
Since $(0,0,0)$ is a smooth point of $\cH$, we necessarily have:
\[r(0,0,0)=\rank \left( \begin{array}{cc}
     g_x(0,0) & g_y(0,0) \end{array} \right) = 1.\] 
On the other hand, $(x,y,0)$ is also a smooth point of $\cH$ for any
$x,y$ small enough, so
\[r(x,y,0)=\rank \left( \begin{array}{cc}
     6x & 0 \\ 0 & 0 \\ g_x(x,y) & g_y(x,y) \end{array} \right) = 1.\] 
 Then $g_y(x,y)=0$ for $x\not= 0$ and since $g$ is analytic in $x$ and
 $y$, this shows that $g$ depends only on the variable $x$ and that $ g_x(0) \not= 0$. Thus
\bea 
\Phi_x(x,y,t)=3x^2+tg_x(x)+t^2h_x(x,y,t)=3x^2+tA(x,y,t),\eea
where $A(x,y,t)$ is analytic and invertible in a neighborhood of $(0,0,0)$.
\bea\Phi_y(x,y,t)=t^2h_y(x,y,t).\eea
 After an analytic change of variables we can rewrite the local equations of $\cH$ on $U$ as $(x'^2-t=0, t^2f(x',y,t)=0)$. As this should be a surface around $(x,y,0)$,  $f \in (x'^2-t)$ and this concludes the proof.

\end{proof}

  By combining the local information on the relative Hilbert scheme $\cH$ and that on the associated normal function derived from the above propositions, we obtain the following

\begin{corollary}
  The general fiber $\cH_t$ of the family $\cH$ is reducible.
\end{corollary}

\begin{proof}

 As before, let $l_1$ denote a point on $\cH_0$ which is one of  the 375 basepoints of $\cH$ and let $l_2$ denote a general point of $\cH_0$. We define the analytic function
$$F: \cH \to \CC$$
$$ (l,t) \to \int_{(l_1, t)}^{(l,t)} \omega(t) $$
where $\omega(t)$ is a holomorphic 3-form on $X_t$, varying holomorphically with $t$. From the previous proposition we know that, in appropriate analytic coordinates $(x,y,t)$, F is given in a neighborhood of $(l_2,0)$ by
$$ x^3 - t x | _{3x^2-t =0}.$$
 Therefore 
$$ F(x,y,t) = x^3-3x^2x = -2x^3 $$
is not constant in  a neighborhood of $(l_2,0)$, because of the branching given by $$x= \left( \frac{t}{3} \right)^{1/2}. $$
 Since $F$ should be constant on each component of $\cH_t$, it follows that $\cH_t$ has at least two different connected components. 

 A more detailed algebraic analysis of the Stein factorization of 
$\cH \to \CC$ will be given in sections 6 and 7.
\end{proof}

\bigskip

\section{The van Geemen lines}

 In \cite{ak2} it is proven that the generic Hilbert scheme $ \cH_{t}$ has a component of dimension 1. The proof is based on identifying a particular set of 5000 points on each  $\cH_{t}$  --the Van Geemen lines-- and proving that, for $t$ generic, the dimension of the tangent space to  $\cH_{t}$ at those points is 1. We know that the generic quintic threefold contains 2875 lines. Since no smooth quintic admits a 2-parameter family of lines, this shows that $\cH_{t}$ has a component of dimension 1. 

  Following is a brief description of how the van Geemen lines occur in the quintic $X_t$. After that,  $\cN_{l_3\backslash X_t}$ and $\cN_{l_3\backslash \cX}$ are computed by the methods introduced in section 1, to the effect that all van Geemen lines are smooth points of $\cH$. Also, all but those over the fibers at the fifth roots of $\frac{2^7}{3}$ are smooth points of their fibers $\cH_t$. 

 Consider lines passing through points of the form: $(1:1:1:a:b)$ and $(1:\xi:\xi^2:0:0)$ in $\PP^4$, where $\xi$ is a complex number satisfying $1+\xi+\xi^2=0$. The conditions for such a line to be contained in one of the $X_t$-s are: \bean  b^5+a^5=27 \mbox{ and } tba=6.  \eean 
  Let $l_3$ denote such a line. There are 10 distinct solutions of the above equation for each $t\not= 0$. Using the action of the group $ G\times S_5$ one obtains a set of 5,000 lines in each $X_t$ (see \cite{ak2}). These will be called the van Geemen lines.

 Consider the Pl\"{u}cker embedding of the Grassmannian $G(2,5)$ in $\PP^9$.
 Notice that all the van Geemen lines correspond to points inside the 10-degree hypersurface in $\PP^9$: $$P=\bigcup_{i,j=0}^4\{p_{ij}=0\}.$$ On the other hand, we know that $\cH_{0} \subset G(2,5)$ consists of 50 quintic curves with multiplicity 2, thus the intersection cycle with P has number 5000. We will see shortly that, indeed, the van Geemen points in $\cH_t$ give exactly the above mentioned hyperplane sections on $\cH_t$.
 
  We next compute the normal bundles to lines $l_3$ in $X_t$, $\cX$, using Lemma 1.1.

  Let $l_3$ be given by the morphism  $ \phi_{3}: \PP^1 \to \PP^4, $ \[\phi_{3}((\alpha : \beta )):= (\alpha + \beta \xi^2:\alpha + \beta \xi :\alpha + \beta :\alpha a:\alpha b), \] with $b^5+a^5=27, tba=6 $.

 As before, one works with the exact sequence:
\bea \begin{CD}
0 @>>> H^0(\cN_{l_3 \backslash X_t}) @>>> \bigoplus_{i=1}^3 H^0(\cO_{l_3}(1))  @>{\psi_{3t}}>> H^0(\cO_{l_3}(5)) 
\end{CD}. \eea 
 Then as before
 \bea \psi_{3t}((a_2,b_2),(a_3,b_3),(a_4,b_4))=\sum_{i=2}^{4}(a_iu+b_iv)\left(\phi_3^*(\frac{\partial F_{t}}{\partial z_i})\right)= \eea
\bea =5(a_2u+b_2v)(-5u^4+10u^3v+4uv^3+v^4)+\eea
\bea +5(a_3u+b_3v)((a^4-tb)u^4-tbuv^3)+5(a_4u+b_4v)((b^4-ta)u^4-tauv^3). 
\eea
 In terms of the basis $(\alpha ,\beta )$ for $H^0(\cO_{l_3}(1))$ and $(\alpha ^5, \alpha^4\beta , ...,\alpha\beta^4, \beta^5)$ for $ H^0(\cO_{l_3}(5))$: 
 \bea \psi_{3t}  = 5 \left( \begin{array}{cccccc}
      -5 & 10 & 0&4&1&0 \\
      0 & -5 & 10&0&4&1 \\
      a^4-tb&0&0&-tb&0&0\\
      0&a^4-tb&0&0&-tb&0\\
      b^4-ta&0&0&-ta&0&0\\
      0&b^4-ta&0&0&-ta&0 \end{array} \right).\eea  
 Notice that since the 3-rd and the 6-th columns of the above matrix are linearly dependent, 
\bea h^0(\cN_{l_3\backslash X_t})=\dim\Ker \psi_{3t} \geq 1 \eea
for all $t\not= 0$. Also notice that if $$ a^5 = b^5 $$ then rows 3 and 5 are linearly dependent and the same is true of rows 4 and 6. In this case one can immediately check that the rank of the above matrix is 4 and so
\bean h^0(\cN_{l_3\backslash X_t})=\dim\Ker \psi_{3t}=2. \eean
 Under the condition $ a^5 = b^5 $, the equations for the van Geemen lines:
 $$b^5+a^5=27, tba=6 $$ yield 
 solutions \[ a= \mu_1\sqrt[5]{\frac{27}{2}},   b= \mu_2\sqrt[5]{\frac{27}{2}}  \] and $$ t = \sqrt[5]{\frac{2^7}{3}}\mu $$
where $\mu_1\mu_2\mu = 1$ and  $\mu_1, \mu_2$ are 5-th roots of unity.

 Otherwise, we check that the determinant
 \bea \left| \begin{array}{ccccc}
      0 & -5 &0&4&1 \\
      a^4-tb&0 &-tb&0&0\\
      0&a^4-tb&0&-tb&0\\
      b^4-ta&0&-ta&0&0\\
      0&b^4-ta&0&-ta&0 \end{array} \right| =
 \left| \begin{array}{ccccc}
      0 & -9 &0&4&1 \\
      a^4&0 &-tb&0&0\\
      0&a^4&0&-tb&0\\
      b^4&0&-ta&0&0\\
      0&b^4&0&-ta&0 \end{array} \right| = \eea
\bea = t^2 (a^{10} - b^{10}) = 27 t^2 (a^5 - b^5) \not=0.\eea  
 Thus for $t \not= \sqrt[5]{\frac{2^7}{3}}\mu, $ 
\bean h^0(\cN_{l_3\backslash X_t})=\dim\Ker \psi_{3t} = 1. \eean

 Similarly, $h^0(\cN_{l_3\backslash \cX})$ can be obtained from the exact sequence
 \be \begin{CD}
H^0(\cN_{l_3\backslash\cX})@>>>H^0(\cN_{l_3\backslash(\PP^4\times\CC)})@>>> H^0(\cN_{\cX \backslash(\PP^4\times\CC})|_{l_3})\\
      @VV{=}V  @VV{=}V  @VV{=}V\\
H^0(\cN_{l_3\backslash\cX})@>>>\bigoplus_{i=1}^3 H^0(\cO_{l_3}(1))\oplus H^0(\cO_{l_3})@>{\tilde{\psi_{3t}}}>>H^0(\cO_{l_3}(5)) 
\end{CD} \ee
where the left horizontal morphisms are injective and $\tilde{\psi_{3t}}$ is obtained from $\psi_{3t}$ by adding the row:
\[ (\begin{array}{cccccc} -ab&0&0&-ab&0&0 \end{array}) \] 
 given by $\phi^*_3(\frac{\partial F_{t}}{\partial t})$, so:
\bean h^0(\cN_{l_3\backslash \cX})=\dim\Ker \tilde{\psi_{3t}}=2. \eean

 These results can be summarized in the following:
 \begin{lemma}
With the above notations, we have:
 \[ \cN_{l_3\backslash X_t}=\cO_{\PP^1}(1)\oplus\cO_{\PP^1}(-3)  \mbox{ for } t=2\sqrt[5]{\frac{4}{3}}\mu, \]
\[ \cN_{l_3\backslash X_t}=\cO_{\PP^1}\oplus\cO_{\PP^1}(-2)  \mbox{ for } t\not=2\sqrt[5]{\frac{4}{3}}\mu, \]
\[\cN_{l_3\backslash \cX}=\cO_{\PP^1}\oplus\cO_{\PP^1}\oplus\cO_{\PP^1}(-2) 
 \mbox{ for all } t.\]

 Thus the relative Hilbert scheme $ \cH$ of the Dwork pencil is a smooth surface at the Van Geemen points; all the fibers $\cH_t$ are smooth at these points except those where t is a fifth root of $\frac{2^7}{3}$.
 \end{lemma}

\bigskip

\section{Algebraic description of the universal line $\cU$ away from $\cU_0$ and the van Geemen lines}
 
\begin{notation}
 As before, let $G(2,5)$ denote the Grassmannian of lines in $\PP^4$. Let
 $F(1,2,5)\subset  G(2,5)\times \PP^4$ denote the flag variety
\bea F(1,2,5) = \{ (\{l\}, y) / \{l\}\in G(2,5), y\in l. \} \eea 

 Let $\cU \subset \cH \times \PP^4$ denote the universal line of $\cH$. 
 Let $\bar{\cU}\subset \bar{\cH} \times \PP^4$ denote the closure of  $\cU$ in $F(1,2,5)$.

 Let $U$ be the projection of $\cU \subset F(1,2,5)\times\PP^1$ to $F(1,2,5)$.
\end{notation} 

\begin{notation}
  For each $i \in \{0,...,4\}$, let $\Psi^i$ denote the rational map from $F(1,2,5)$ to $\PP^3$, defined by
\bea \Psi^i((\{l\}, y))= \left(\frac{x_0}{y_0}:...:\hat{\frac{x_i}{y_i}}:...:\frac{x_4}{y_4} \right), \eea
 where $y=(y_0:...:y_4)\in l \subset \PP^4$ and $x=(x_0:...:x_4)\in l\cap H_i$.
\end{notation}

 The restriction on $U$ of the rational map $\Psi^i$ will be the main tool in our algebraic study of $\cU$ away from the Van Geemen lines and from the lines in $X_0$. We will denote this restriction by $\psi^i$.

 We can understand the map $\Psi^i$ geometrically as follows: for a
 given line $l \subset \PP^4 \backslash
 \left(\bigcup_{j=0}^4H_j\right)$ and a fixed point $y\in l$, we
 consider $y$ to be the point at $\infty$ on $l$, and $l\cap H_i$ to
 be the $0$ point on $l$. Modulo the action of $\CC^*$, this gives a
 parametrization on $l$. Thus, 4-tuples of points on $l$, none of
 which is $\infty$, will correspond via this structure on $l$ to
 points of $\PP^3$. We let $ \Psi^i((\{l\}, y))\in \PP^3$ be given by
 the 4-tuple of intersection points $\{ l\cap H_j \}$, where $j\not=
 i$. Note that this construction only works for  $l \not\subset
 \bigcup_{j=0}^4H_j$ and $y \in \PP^4 \backslash
 \left(\bigcup_{j=0}^4H_j\right)$. It extends algebraically to $ (\{l\}, y)$ with  $y\in H_j\backslash  \left(\bigcup_{k\not= j}H_k\right)$, as long as  $l \not\subset \bigcup_{j\not=i}H_j $.

 We may think of $F(1,2,5)$ as the projectivization of the tangent bundle $T_{\PP^4}$. Over $$(\CC^*)^4 =\PP^4 \backslash \left(\bigcup_{j=0}^4H_j\right),$$ this bundle is trivial and any of the maps $\Psi^i$ gives a trivialization of it.

 \begin{notation} Let $\cU' \subset \cU$ be defined as:
\bea \cU' =  \{ (\{l\}, y) \in \cU /  \prod_{j=0}^4 y_j \not=0 \}. \eea
 Let $U' \subset U$ be the projection of $\cU'$ on $F(1,2,5)$.
 \end{notation}

  \begin{notation}  Consider the Pl\"{u}cker embedding of $G(2,5)$ in $\PP^9$, with coordinates $p_{ij}$.  Let $\cH' \subset \cH$ be defined as the preimage in $\cH$ of 
$$ G(2,5) \backslash \left( \bigcup_{i,j=0}^4 \{p_{ij} = 0 \} \right).$$
 \end{notation}

 The following notations will prove very useful in better understanding the subset $\cU' \subset \cU$ and the morphisms 
\bea \psi^i |_{U'} : U' \to \PP^3. \eea

\begin{notation}

 Let $(\{l\}, y) \in F(1,2,5)$ be such that $\prod_{j=0}^4 y_j \not=0$. For each $x=(x_0:...:x_4) \in l$, let
\bea  u_j = \frac{x_j}{y_j},  \forall j \in\{0,...,4\}. \eea
 
 Following the discussion above, if $x$ varies in one of the coordinate hyperplanes $H_i \subset\PP^4$, then $((u_0:...:\hat{u_i}:...:u_4), (y_0:...:y_4)) \in \PP^3\times\PP^4$ give coordinates for the restriction of the $\PP^3$-bundle $F(1,2,5)$ over $\PP^4$ to   $(\CC^*)^4 =\PP^4 \backslash \left(\bigcup_{j=0}^4H_j\right)$. 
 Notice also that for $x\not= y$ varying in a hyperplane $H_i\subset\PP^4$, the projective coordinates $(u_0:...:\hat{u_i}:...:u_4)$  together with $(y_0^5:...:y_4^5)$ exactly determine the orbit of the line $l$ by the action of the group $\hat{G}$.
Let
 \[s_0(u)=1, s_1(u)=\sum_{i=0}^{4}u_i,  s_2(u)=\sum_{i\not= j; i,j=0}^4u_iu_j,  \mbox{ etc. }\]
 be the fundamental symmetrical polynomials in $(u_0,...,u_4)$. 

For any $k\in\{0,...4\}$, let
\[\sigma_k(u, y)=\sum_{i=0}^{4} u_i^k y_i^5.\]
 For $i, j\in\{0,...4\}$, we let $$u^i = (u_0:...:0:...:u_4)$$
where $u_i=0$ and $$u^{ij} = (u_0:...:0:...:0:...:u_4)$$
where $u_i=u_j=0$. 
 
 Also, let
 $$\delta(u):=\prod_{j>k} (u_j-u_k),$$
 $$\delta(u^i):=\prod_{j,k\not= i;  j>k} (u_j-u_k),$$
etc.
\end{notation}

\begin{notation}  If $x\not= y$, we may consider the line $l$ passing through $x$ and $y$ to be parametrized with coordinates $(\alpha :\beta )$ and embedded in $\PP^4$ by the morphism $\phi : \PP^1 \to \PP^4 $,
\[\phi((\alpha :\beta )):= (x_0\alpha +y_0\beta :x_1\alpha +y_1\beta :x_2\alpha +y_2\beta :x_3\alpha +y_3\beta :x_4\alpha +y_4 \beta ).\]
 In the new coordinates $((y_j)_j, (u_j)_j)$: 
\bean \phi((\alpha :\beta )) = (y_0(u_0\alpha +\beta ):y_1(u_1\alpha +\beta ):...:y_4(u_4\alpha +\beta )). \eean\end{notation}

\begin{notation}
Let $M(u)$ denote the matrix 
\bea M(u)= \left( \begin{array}{ccccc} 1&1&1&1&1 \\
                          u_0 & u_1 & u_2 & u_3 & u_4 \\
                          u_0^2 & u_1^2 & u_2^2 & u_3^2 & u_4^2 \\
                          u_0^3 & u_1^3 & u_2^3 & u_3^3 & u_4^3 \\
                          u_0^4 & u_1^4 & u_2^4 & u_3^4 & u_4^4 \end{array} \right). \eea
 Let $$ Y(y)=  \left( \begin{array}{c} y_0^5 \\y_1^5\\y_2^5 \\ y_3^5 \\ y_4^5 \end{array} \right)  \mbox{ and }  C(u)= \left( \begin{array}{c} s_0(u) \\ \frac{1}{5}s_1(u) \\ \frac{1}{10}s_2(u) \\ \frac{1}{10}s_3(u) \\ \frac{1}{5}s_4(u)\end{array} \right). $$
\end{notation}

\begin{proposition}
  $\cU'$ is a 3-dimensional quasi-projective variety, reduced and complete intersection in $F(1,2,5)_{| (\CC^*)^4 }$, given, via the trivialization $\Psi_i$ of $F(1,2,5)_{| (\CC^*)^4 }$, by the following matrix equation on  $ \PP^3\times(\CC^*)^4 $:
\bean M(u^i)\cdot Y(y) = 5t\prod_{j=0}^4 y_j  C(u^i). \eean 
\end{proposition}

\begin{proof}
  Consider the pair 
$$(\{l\}, y) \in F(1,2,5)|_{(\CC^*)^4},$$ with the line $l$
parametrized and embedded in $\PP^4$ by the morphism $\phi : \PP^1 \to
\PP^4 $ given in equation (4.1). The line $l$ is in one of the
$X_t$-s, for some $t\in\CC$, if $\phi^*(F_t ) \equiv 0 $
as a degree 5 homogeneous polynomial in $(\alpha :\beta )$.
 This gives six equations in $u=(u_0:...:u_4)$, $y=(y_0:..:y_4)$:
\bean \left\{ \begin{array}{l}
 \left(\begin{array}{l}5\\k\end{array}\right)\sigma_k(u,y)-5ts_k\prod_{j=0}^4 y_j =0,  \mbox{ for } k\in \{ 0,..,5\}.
 \end{array} \right. \eean 
Equivalently, in matrix form:
\bean \left( \begin{array}{ccccc} 1&1&1&1&1 \\
                          u_0 & u_1 & u_2 & u_3 & u_4 \\
                          u_0^2 & u_1^2 & u_2^2 & u_3^2 & u_4^2 \\
                          u_0^3 & u_1^3 & u_2^3 & u_3^3 & u_4^3 \\
                          u_0^4 & u_1^4 & u_2^4 & u_3^4 & u_4^4 \\
                          u_0^5 & u_1^5 & u_2^5 & u_3^5 & u_4^5 
 \end{array} \right)  \left( \begin{array}{c} y_0^5 \\y_1^5\\y_2^5 \\ y_3^5 \\ y_4^5 \end{array} \right) =
 5t\prod_{j=0}^4 y_j \left( \begin{array}{c} s_0(u) \\ \frac{1}{5}s_1(u) \\ \frac{1}{10}s_2(u) \\ \frac{1}{10}s_3(u) \\ \frac{1}{5}s_4(u) \\ s_5(u)\end{array} \right). \eean 
 As long as $5t\prod_{j=0}^4 y_j \not= 0$, one can forget this factor since $(y_0:...:y_4)$ are given modulo the action of $\CC^*$.

 Notice that the sixth row $r_5$ of the augmented matrix of the system:
\bea \left( \begin{array}{cccccc}  u_0^5 & u_1^5 & u_2^5 & u_3^5 & u_4^5 & s_5(u) \end{array} \right) \eea
is a linear combination of the first five rows $(r_i)_{i\in\{0,...,4\}}$:
\bea r_5 - s_1(u) r_4 + s_2(u) r_3 - s_3(u) r_2 + s_4(u) r_1 -s_5(u) r_0=0. \eea 
 Thus after setting $u_i=0$, one gets 5 equations defining a reduced complete intersection in $\PP^3\times (\CC^*)^4$.
 \end{proof}

\begin{corollary}
 Let $y\in X_t$ for some $t\in\CC$ satisfy  $$\prod_{j=0}^4 y_j \not=0.$$ Then a line $l\subset \PP^4$ has contact of order $\geq 5$ with $X_t$ at the point $y$ if and only if 
\bea l \subset X_t \eea
\end{corollary}

\begin{proof}
  Looking back at the way equations (4.4) were obtained, one notices that the first $k$ equations in (4.4) define the subset
$$\{(\{l\},y) \in F(1,2,5)|_{(\CC^*)^4} / \mbox{ l has contact of order}\geq k \mbox{ with } X_t \mbox{ at } y \}$$ 
of $F(1,2,5)|_{(\CC^*)^4}$. Thus the corollary rephrases the statement that the 6-th equation in (4.4) is a combination of the other five.
\end{proof}

\begin{theorem}
 The inverse image of  $\cH' / \hat{G}$ in $ \cU / \hat{G}$ is a smooth, irreducible, rational threefold. Thus $\cH' / \hat{G}$ is a smooth
irreducible unirational, therefore rational, surface.

\end{theorem}

\begin{proof} 
 If $u$ is an element of $\PP^4 \backslash \left(\bigcup_{j\not=k} \{u_j = u_k\} \right)$, then the matrix $M(u)$ is invertible.
 In this case, equation (4.2) can be rewritten as
 \bean Y(y) =5t\prod_{j=0}^4 y_j  M(u^i)^{-1} C(u^i).  \eean
 Recall that $Y(u)$ is given only in terms of $y^5:=(y_0^5:...:y_4^5)$. As long as $t\prod_{j=0}^4 y_j \not= 0$, one may disregard the $t\prod_{j=0}^4 y_j$. As a consequence, equation (4.5) defines a birational transformation
 $$\gamma^i : \PP^3 \to  \cU/\hat{G},$$
$$ \gamma^i ( u ) = ( y^5, u ). $$

 If one thinks of $G(2,5)$ as embedded in $\PP^9$ by the Pl\"{u}cker embedding, then, in terms of the coordinates $((y_j)_j,(u_j)_j)$, $$ p_{ij} = y_iy_j(u_i-u_j). $$ So, if $y$ is such that $\prod_{i=0}^4 y_i \not= 0$, then the conditions $y\in l$, $\{l\} \in \cH'$ are equivalent to the fact that $M(u)$ is invertible. The map $\gamma^i$ defined above is thus a smooth morphism at such pairs $(l,\{y\})$ and therefore the images of these points are smooth in $\cU/\hat{G}$. The points of $\cU$ standing over these points of $\cU/\hat{G}$ are also smooth.
\end{proof}

 Next we will take a closer look at the lines excluded by the previous corollary:

\begin{theorem}
 Let $l\subset X_t$ be a line such that $\{l\} \in \cH \backslash \cH'$. Then either $l\subset X_0$ or $l$ is a Van Geemen line.
\end{theorem}

\begin{proof}
 For symmetry reasons, it is enough to check the proposition for $l$ with $p_{01} = 0$. Assume that $l \not\subset X_0$. Then one can choose a point $y$ on $l$ such that $\prod_{i=0}^4 y_i \not= 0$. Indeed, one can always do that as long as $l$ is not included in one of the coordinate hyperplanes $H_i$. But if $l\subset H_i$ and  $l\subset X_t$ for some $t\in\CC$, then $l$ is in the base locus of the pencil $\cX$ and thus also in $X_0$. 

 In the proof of Corollary 4.4 it is shown that the condition $p_{01} = 0$ becomes $u_0 = u_1 $. By considering the map $\psi^0$, one may assume $$u_0 = u_1 =0.$$ The proposition can be proven by analyzing equation (4.2) in this context.

 Because $u_0 = u_1 =0$, all the other $u_j$-s satisfy the following:
$$ u_j^{k+3} - s_1(u^{01}) u_j^{k+2} + s_2(u^{01}) u_j^{k+1} - s_3(u^{01}) u_j^{k} =0. $$
 In particular, for $k=1$ one obtains that the last 4 rows of the matrix $M(u)=M(u^{01})$ are linearly dependent and thus, the last 4 rows of $C(u)$ should necessarily satisfy the same dependence relation:
 $$\frac{s_4(u^{01})}{5} -s_1(u^{01}) \frac{s_3(u^{01})}{10} +s_2(u^{01}) \frac{s_2(u^{01})}{10} -s_3(u^{01}) \frac{s_1(u^{01})}{5} = 0.$$
 Because $s_4(u^{01}) = 0$, this relation simplifies to:
\bean s_2^2(u^{01})-3s_1(u^{01}) s_3(u^{01}) = 0,\eean
\bean ( u_3u_4 + \xi u_2u_4 + \xi^2 u_2u_3 )( u_3u_4 + \xi^2 u_2u_4 + \xi u_2u_3 ) = 0 \eean 
 where $\xi$ is a primitive third root of unity.

 For symmetry reasons it is enough to study the case 
\bean u_3u_4 + \xi u_2u_4 + \xi^2 u_2u_3 =0.\eean 
 If two of the $u_i$-s are equal, then either all three are equal or those two are actually 0. If all three are equal, then $l$ must be passing through a point $(\mu_0: -\mu_1 : 0: 0: 0)$ with $\mu_0^5 = \mu_1^5 = 1$, and also through a point $(0: 0: x_2:x_3:x_4)$, because of the assumption $p_{01}=0$. In this case $l \subset X_0$.

 Notice also that if any one of the coordinates $(u_2,u_3,u_4)$ is 0, then yet another one should be 0. But it is easy to check that there is no line inside $\cX$ passing through a point four of whose coordinates are 0. Thus one obtains that $u_2u_3u_4 \not= 0$. This, together with equation (4.8), justifies a change of coordinates to $v^{01}:=(v_2:v_3:v_4)$, where $$v_i = \frac{1}{u_i}. $$ Equation (4.8) becomes:
 \bean v_2 + \xi v_3 + \xi^2 v_4 =0. \eean  

 The $u$-terms involved in the algebraic description (4.5) of $\cU'$ transform as follows under the change $v_i = \frac{1}{u_i}$:
\bean s_k (u') = \frac{s_{top - k}(v')}{s_{top}(v')}, \eean
where $top = 5$ for $u'=u$ and $v'=v$, $top = 4$ for $u'=u^i$ and $v'=v^i$ etc.
\bean \delta (u') =(-1)^{\frac{ top ( top - 1)}{2}} \frac{\delta (v')}{s_{top}(v')^{top -1}}. \eean

 Given $(v_2:v_3:v_4)$, equations (4.5) well determine  $(y_2^5:y_3^5:y_4^5)$ (up to a common factor). Indeed, after dividing by the $t\prod_{i=0}^4 y_i$- term:
$$ Y^{01}(y^{01}) = M^{01}(u^{01})^{-1} C^{01}(u^{01}), $$
where  $$ Y^{01}(y^{01}) := \left(\begin{array}{c} y_2^5 \\y_3^5 \\ y_4^5 \end{array} \right),  C^{01}(u^{01}) := \left(\begin{array}{c} s_1(u^{01}) \\ \frac{1}{2}s_2(u^{01}) \\ \frac{1}{2}s_3(u^{01})    \end{array} \right), $$
$$ M^{01}(u^{01}) := \left( \begin{array}{ccc}
                         u_2 & u_3 & u_4 \\
                          u_2^2 & u_3^2 & u_4^2 \\
                          u_2^3 & u_3^3 & u_4^3 
 \end{array} \right).  $$
 Alternatively, following formulas (4.10), and (4.11):
\bean Y^{01}(y^{01}) = M^{01}(v^{01})^{-1} C^{01}(v^{01}) \eean 
with $$C^{01}(v^{01}) = \frac{1}{s_3(v^{01})}\left(\begin{array}{c} s_2(v^{01}) \\ \frac{1}{2}s_1(v^{01}) \\ \frac{1}{2}s_0(v^{01}) \end{array} \right) $$  and the product $$\delta (v^{01})M^{01}(v^{01})^{-1}$$ equal to $$\left( \begin{array}{ccc} v_2^3\delta(v^{012})s_0(v^{012}) &
 -v_2^3 \delta(v^{012})s_1(v^{012}) & v_2^3 \delta(v^{012})s_2(v^{012}) \\
 -v_3^3 \delta(v^{013})s_0(v^{013}) &
 + v_3^3 \delta(v^{013})s_1(v^{013}) &  -v_3^3 \delta(v^{013})s_2(v^{013}) \\
 v_4^3 \delta(v^{014})s_0(v^{014}) &
 -v_4^3 \delta(v^{014})s_1(v^{014}) & v_4^3 \delta(v^{014})s_2(v^{014}) 
 \end{array} \right). $$ 
 Thus one derives:
\bea y_j^5 =(-1)^j \frac{\delta (v^{01j}) v_j^3}{\delta (v^{01}) s_3(v^{01})}\left( s_2(v^{01})-\frac{1}{2}s_1(v^{01j})s_1(v^{01})+\frac{1}{2} s_2(v^{01j})\right) \eea for $j \in \{ 2, 3, 4 \}$.
 Using the identities:
$$ s_1(v^{01j}) = s_1(v^{01}) - v_j, $$
$$ s_2(v^{01j}) =  s_2(v^{01}) - v_js_1(v^{01}) + v_j^2, $$
 one can further simplify to  
\bea y_j^5 = \frac{(-1)^j\delta (v^{01j}) v_j^3}{2 \delta (v^{01} s_3(v^{01})}\left( 3s_2(v^{01}) - s_1^2(v^{01}) + v_j^2 \right). \eea
 But formula (4.9) implies $ 3s_2(v^{01}) - s_1^2(v^{01}) =0$ and then
\bea y_j^5 = \frac{(-1)^j\delta (v^{01j}) v_j^5}{2 \delta (v^{01}) s_3(v^{01})}, \eea
or equivalently,
\bea x_j^5 = \frac{y_j^5}{ v_j^5} = \frac{(-1)^j\delta (v^{01j})}{2 \delta (v^{01}) s_3(v^{01})} \eea
 for  $j \in \{ 2, 3, 4 \}$. Thus one obtains:
\bea (x_0^5:x_1^5:x_2^5:x_3^5:x_4^5) = (0 :0 :v_4-v_3: v_2-v_4: v_3-v_2)= \eea 
                               \bea      = (0 :0 :1: \xi : \xi^2), \eea 
since $v^{01}\not= (1:1:1) $ lies in the plane $v_2+\xi v_3 +\xi^2 v_4=0$. Then
\bean (x_0:x_1:x_2:x_3:x_4) = (0:0: \mu_2: \xi^2\mu_3: \xi\mu_4) \eean
for some fifth roots of unity $\mu_2, \mu_3, \mu_4 $. From here,
\bean (y_2:y_3:y_4) = (v_2x_2:v_3x_3:v_4x_4) = (v_2\mu_2: v_3\xi^2\mu_3: v_4\xi\mu_4). \eean
 In view of formula (4.9), the last statement leads to the following: If $l$ satisfies $p_{01}=0$, then its projection on the last three coordinates:
\bea \PP^1 &\rightarrow & \PP^2, \\
      (\alpha :\beta ) &\rightarrow &(\phi_2(\alpha ,\beta ):\phi_3(\alpha ,\beta ):\phi_4(\alpha ,\beta ) ), \eea
is a line of equation 
\bea \mu_2^4 z_2 +\xi^2 \mu_3^4 z_3 +\xi \mu_3^4 z_3 = 0 \eea
for some fifth roots of unity $\mu_2, \mu_3, \mu_4,$ or $$\mu_2^4 z_2 +\xi \mu_3^4 z_3 +\xi^2 \mu_3^4 z_3 = 0, $$ if one looks at the other component of the locus of such lines pointed out in formula (4.7).
 In particular, such a line contains a point with $$(y_2:y_3:y_4)=(\mu_2: \mu_3: \mu_4).$$ This, together with formula (4.13), are enough to characterize a Van Geemen line.
\end{proof}

 Putting together the information on smoothness gathered along these 4 sections, one may conclude:

\begin{proposition}
 Outside the lines in the base locus of the family $\cX$, $\cU$ is a smooth threefold.
\end{proposition}

\section{The surface component $\cS$ of $\cH$}

Putting together the results of Sections 1 through 4, one can prove: 

\begin{theorem}
The irreducible components of the relative Hilbert scheme $\cH$ are given by:
\[ \cH = \cS \cup \bigcup_{i=1}^{375} Z_i \]
where $\cS$ is a smooth quasi-projective surface proper over $\CC$ and for each $i$,\\
$\pi_{| Z_i} :Z_i \to \CC$ is an isomorphism.
 The Hilbert polynomial of the family $\cS$ inside G(2,5) is
\[ p_t(n)=500n-1250. \]
\end{theorem}

\begin{proof}
 Consider the Pl\"{u}cker embedding $G(2,5) \subset \PP^9$. Let  $$ T
 = \PP^9 \backslash \left( \bigcup_{i,j=0}^4 \{p_{ij} = 0 \}\right)$$ denote the torus in $\PP^9$:
 
 Theorem 4.4 shows that for each $t\not= 0$, the intersection of $\cH_t$ with the border of the torus T consists of the 5000 points (2500 for $t=\ts$), representing the van Geemen lines of the quintic, together with 375 isolated points. 
 Following the description in Section 3, the union of all the Van Geemen points of the relative Hilbert scheme $\cH$ gives a set of 500 smooth irreducible curves in $ \PP^2 \times \CC$ satisfying equations like these:
\[ \left\{ \begin{array}{l}
b^5+a^5-27c^5=0\\
tba-6c^2=0 
\end{array} \right., \] 
 where $(a:b:c)$ are homogeneous coordinates in $\PP^2$. Each of these curves is a 10-1 cover of  $\CC$, branched over $0$ and $\ts $. From the Propositions 2.2 and 2.3 one gathers that there is only one irreducible 2-dimensional component $\cS$ of $\cH$ containing $\cH_0$ and moreover, this is the only irreducible 2-dimensional component of $\cH$ which intersects $\cH_0$. But since all the irreducible curves of van Geemen points intersect $\cH_0$, it follows that they are all contained in $\cS$.

 On the other hand, as a result Section 4, the preimage $\cH'$ of
 $G(2,5)\cap T$ in $\cH$ is also irreducible. Indeed, Theorem 4.3 states that $\cH'/\hat{G}$ is irreducible; a simple exercise starting from equations (4.5) shows that $\cH'$ itself is irreducible. One may also find a further verification of this fact in formula (6.6), which explicitly describes the map $\gamma^i$. 

 Theorem 4.4 shows that $\cH\backslash\cH'$ consists of the lines in the base locus of $\cX$, plus those in $\cH_0$ and the van Geemen lines, contained, as seen above, in $\cS$. Form the structure of $\cH$:
\[ \cH = \cS \cup \bigcup_{i=1}^{375} Z_i \]
where $\cS$ is a smooth quasi-projective surface proper over $\CC$ and for each $i$, $Z_i$ corresponds to some line in the base locus of $\cX$, so
$\pi_{| Z_i} :Z_i \to \CC$ is an isomorphism.

The Hilbert polynomial of the family $(\cS_t)_t$ can be computed as
\bea p(n)=2p_{\cH_{0\mbox{ red }}}(n)=2(50(5n-5)-375)=500n-1250 \eea
since each $C_{i,j,\mu}$ is embedded in G(2,5) as a degree 5 plane curve, each of the 375 base-locus points is at the intersection of exactly 2 of the components of $\cH_{0\mbox{ red }}$, with   intersection multiplicity 1, while the whole  $\cH_{0\mbox{ red }}$ is embedded in $\cH$ with embedding multiplicity 2.
\end{proof}

\begin{theorem}
 The quotient $\cS/\hat{G}$ of the surface $\cS$ is a smooth, irreducible, quasi-projective rational surface, proper over $\CC$.
\end{theorem}

\begin{proof}
 The finite group $\hat{G}$ acts on the smooth surface $\cS$. One needs to check whether the quotients of those points having non-trivial isotropy groups under this action are smooth points of  $\cS/\hat{G}$. First recall from Theorem 4.3 that the open subset $\cH'/\hat{G}$ of $\cS/\hat{G}$ is smooth and rational. The complement of $\cH'$ in $\cS$ consists of the van Geemen lines and the lines in $\cH_0$. Each of the van Geemen lines has a trivial isotropy group, as can be easily seen from their definition. Consider now a line $l_1$ in the base locus of the family $\cX$, and a line $l_{2}$ contained in $X_{0}$, but not in the base locus.  The point $\{l_{2}\}\in\cH_{0}$ has order 5 isotropy group $G_{2}$: if, for example, $l_{2}$ passes through $(1:-1:0:0:0)$ and $(0:0:x_2:x_{3}:x_4) \in \PP^{4}$, then $$G_{2}=\{(1:1:\mu:\mu:\mu) / \mu^5=1\}.$$
  The point $\{l_{1}\}\in\cH_{0}$ has order 25 isotropy group $G_{1}$: if, for example, $l_{1}$ passes through $(1:-1:0:0:0)$ and $(0:0:1:-1:0)\in \PP^{4}$, then $$G_{1}=\{(1:1:\mu_{1}:\mu_{1}:\mu_{2}) / \mu_{1}^5=1, \mu_{2}^5=1\},$$
with two generators $(1:1:\mu:\mu:\mu)$ and $(\mu:\mu:1:1:\mu)$. Each
  of these generators gives a pseudoreflection of $\cS$, i.e. the set of points in $\cS$ fixed by it is a divisor of $\cS$. Indeed, the set of points in $\cS$ fixed by $(1:1:\mu:\mu:\mu)$ corresponds to the lines in the cones $C_{0,1,\nu}$, of vertices 
$$(1:-\nu:0:0:0),$$
 $\nu^{5}=1$, over the curve $$\{(0:0:x_2:x_{3}:x_4) /
x_{2}^{5}+x_{3}^{5}+x_{4}^{5}=0\}.$$ Following \cite{st}, the above
property is sufficient for the images of  $l_1$ and $l_{2}$ to be
smooth points of  $\cS/\hat{G}$. Thus  $\cS/\hat{G}$ is smooth everywhere.

\end{proof}

\bigskip

\section{Properties of the fiber $( \cS / \hat{G})_{w}$}

 We are now ready to prove:
\begin{proposition}
 The family $( \cS / \hat{G})$ has double fiber at $w = 0, \frac{2^7}{5}$.

The fiber $( \cS / \hat{G})_{w}$ at $w \not= 0, \frac{2^7}{5}$ consists of two isomorphic connected curves.
\end{proposition}

\begin{proof}
  Recall from Corollary 4.3 that the following rational map
$$  \hat{\psi}^i : \cU/\hat{G}  \rightarrow \PP^3 $$
is birational. Let $U^i_w$ denote the closure of  $\hat{\psi}^i( (\cU/\hat{G})_w)$. Starting from formula (4.5), we will derive the equation of the surface  $U^i_w$ in $\PP^3$ and then will prove that $U^i_w$ has two irreducible surface components.

 First, one can explicitly write down the inverse $\gamma^i$ of the above map. Although in this context the $i$-th coordinate $u_i=0$, we will treat it here as all the other coordinates, for symmetry purposes. If $u_j\not= u_k$ for all $j,k \in\{0,...,4\}$, the matrix $M^{-1}(u)$ has elements
$$  M^{-1}(u)_{jk} = (-1)^{j+k} \frac{\delta (u^j)}{\delta (u)} s_{4-k}(u^j) , $$
where $j,k \in \{0,...,4\}$. (Indeed, starting from the identity
$$ \sum_{k=0}^4 (-1)^ks_{4-k}(u^j) u_l^k = \prod_{h\not= j} (u_l-u_h) = \delta_{jl} (-1)^j \frac{\delta (u)}{\delta (u^j)}, $$ one finds that
 $$ \sum_{k=0}^4M^{-1}(u)_{jk} M(u)_{kl} = \sum_{k=0}^4(-1)^{j+k} \frac{\delta (u^j)}{\delta (u)} s_{4-k}(u^j)u_l^k = \delta_{jl}, $$
 where $\delta_{jl}$ is the Kr\"oenecker delta symbol.)

By formula (4.5), $\gamma^i$ is given by 
\bea y_j^5 = \sum_{k=0}^{4}  M^{-1}(u)_{jk} C(u,y)_k, \eea so
 \bea   y_j^5 = (5t\prod_{l=0}^{4}y_l) \frac{\delta (u^j)}{\delta (u)} \sum_{k=0}^{4} \frac{(-1)^{j+k}}{\left(\begin{array}{l}5\\ k\end{array}\right)}  s_{4-k}(u^j) s_k(u). \eea
 Use the identity 
\bea s_k(u) = s_k(u^j) + u_j s_{k-1}(u^j)  \eea
to further simplify:
\bea y_j^5 &=& (5t\prod_{l=0}^{4}y_l) \frac{\delta (u^j)}{\delta (u)} \sum_{k=0}^{4} \frac{(-1)^{j+k}}{\left(\begin{array}{l}5\\ k\end{array}\right)}  s_{4-k}(u^j)\left( s_k(u^j) + u_j s_{k-1}(u^j)\right) = \\
          &=& (5t\prod_{l=0}^{4}y_l) \frac{\delta (u^j)}{\delta (u)}\frac{(-1)^j}{10} \left( s_2^2(u^j)-3s_1(u^j)s_3(u^j)+12s_4(u^j)\right) . \eea
\begin{notation} 
Let $g(u^j):=s_2^2(u^j)-3s_1(u^j)s_3(u^j)+12s_4(u^j)$.
\end{notation}

 The following lemma can be verified by direct computations:

\begin{lemma}
Let $\{h,i,j,k,l\} = \{0,...,4\}$. $g(u^i)$ has the following proprieties:
\bean g(u^i)|_{\{ u_j=u_i\} }=g(u^j)|_{\{ u_j=u_i\} } \eean
\bean g(u^i)|_{\{ u_j=u_k\} }=(u_l-u_j)^2(u_h-u_j)^2 \eean
\bean \frac{\partial g(u^i)}{\partial u_j}|_{\{ u_j=u_k\} }=(u_l-u_j)(u_h-u_j)(u_l+u_h-2u_j) \eean
\bean \frac{\partial^2 g(u^i)}{\partial u_j^2} = \frac{\partial^2 g(u^j)}{\partial u_i^2} =2 \sum_{k\in\{0,..,4\}\backslash\{i,j\}}(u_h-u_k)(u_l-u_k) \eean
\bean \frac{\partial^2 g(u^i)}{\partial u_j^2} |_{\{ u_k=u_l\} }=2(u_h-u_k)^2  \eean
\end{lemma}

 To recap, we have shown that the inverse $\gamma^i$ of $\hat{\psi}^i$ is given by:
\bean  y_j^5 =  5t(\prod_{l=0}^{4}y_l) \frac{\delta (u^j)}{\delta (u)}\frac{(-1)^j}{10} g(u^j). \eean 
 By taking product of the equations (6.6) when $j$ varies from 0 to 4 and then simplifying the $\prod_{j=0}^{4} y_j^5$-factor, one finds the pullback by $\gamma^i$ of the equation for $(\cU/\hat{G})_w$, on the set  $\PP^3 \backslash \left( \bigcup_{i\not=j} \{u_i = u_j \}\right)$:
\bea 1 =\frac{w}{2^5} \prod_{j=0}^{4} \frac{\delta (u^j)}{\delta (u)}g(u^j), \eea
where $w = t^5$. After simplification, here is the equation of the surface $U^i_w$ in $\PP^3$:
\bean \delta^2(u) =\frac{w}{2^5} \prod_{j=0}^{4} g(u^j) \eean

\begin{notation}
Let \bea G(u) := \prod_{j=0}^{4} g(u^j). \eea 
\end{notation}

 To prove the existence of two irreducible components of $U^i_w$, we will find a symmetric, homogeneous degree 10 polynomial $P(u)$ such that 
\bean  \prod_{j=0}^{4} g(u^j) -\frac{3}{4}\delta^2(u)=P^2(u). \eean
 Equations (6.7) and (6.8) will thus show that the fiber at $w=\frac{2^7}{3}$ appears with multiplicity 2 in the family $((\cU/\hat{G})_w)_w$. Also, equation (6.7) may be rewritten as 
\bean \left(1-\frac{3w}{2^7} \right)\delta^2(u)-\frac{w}{2^5}P^2(u)=0, \eean
which factors as the product of two equations for each $w\in \CC$. Moreover, for generic $w$, each of these equations define an irreducible surface in $\PP^3$.

 In order to write $P(u)$ explicitly, one can consider $\frac{G(u)}{\delta^2(u)}$ as a rational function in the variable $u_0$ and decompose it into a sum of simple fractions:
\bean \frac{G(u)}{\delta^2(u)}=\sum_{i=1}^4\frac{A_i(u^0)}{(u_0-u_i)^2}+\sum_{i=1}^4\frac{B_i(u^0)}{(u_0-u_i)}+C(u^0),\eean
with coefficients satisfying:
\bea && A_i(u^0)=\frac{1}{\delta^2(u^0)\prod_{j\not=0,i}(u_i-u_j)^2} G(u)|_{u_0=u_i},\\
     && B_i(u^0)+A_i(u^0)\sum_{j\not=0,i}\frac{2}{u_i-u_j}=\frac{1}{\delta_0^2\prod_{j\not=0,i}(u_i-u_j)^2} \frac{\partial G(u)}{\partial u_0}\left|_{u_0=u_i} \right. , \\ 
     && C(u^0)=\frac{1}{2^4\delta^2(u^0) }g(u^0)\prod_{i=1}^4\frac{\partial^2 g(u^i)}{\partial u_0^2}. \eea
Using properties (6.1) - (6.3) in lemma 6.2, one finds:
\bean  A_i(u^0)=\frac{g^2(u^0)}{\delta^2(u^{0i})} \mbox{ and }
      B_i(u^0)=\frac{1}{\delta_{0i}^2}g(u^0)\frac{\partial g(u^0)}{\partial u_i}.\eean
 Let $a_i(u^0)=\frac{g(u^0)}{\delta (u^{0i})}$. Next one checks that 
\bea  \frac{G(u)}{\delta^2(u)} - \frac{3}{4}
=\sum_{i=1}^4\frac{A_i(u^0)}{(u_0-u_i)^2}+\sum_{i=1}^4\frac{B_i(u^0)}{(u_0-u_i)}+C(u^0)- \frac{3}{4}= \eea 
\bea   =\left(\sum_{i=\bar{1,4}}\frac{a_i(u^0)}{(u_0-u_i)}+E(u^0)\right)^2 \eea where 
\bea E(u^0)=E_i(u^0)=\frac{B_i(u^0)}{2a_i(u^0)}-\sum_{ j\not=0,i}\frac{a_j(u^0)}{(u_i-u_j)} \eea 
 does not  depend on $i$ and 
\bean E^2(u^0)=C(u^0) - \frac{3}{4}. \eean 
 Indeed, one can verify the last two properties by again decomposing into sums of simple fractions in new variables $u_i$. It is enough to do it for $i=1$:
\bea && E_1(u^0)=\frac{1}{2\delta (u^{01})}\frac{\partial g(u^0)}{\partial u_1}- \sum_{j=2}^4\frac{g(u^0)}{\delta (u^{0j})(u_1-u_j)} = \\
&& = \frac{1}{\delta (u^{01})} \left(\frac{1}{2}\frac{\partial g(u^0)}{\partial u_1}-\sum_{j=2}^4\frac{g(u^0)(u_i-u_j)(u_k-u_j)}{(u_1-u_2)(u_1-u_3)(u_1-u_4)} \right) =\\
&& = \frac{1}{\delta (u^{01})} \left(\frac{1}{2}\frac{\partial g(u^0)}{\partial u_1}-\frac{g(u^0)}{2(u_1-u_2)(u_1-u_3)(u_1-u_4)} \frac{\partial^2 g(u^0)}{\partial u_1^2}\right) =\\
&& = \frac{1}{\delta (u^{01})} \left(\frac{1}{2}\frac{\partial g(u^0)}{\partial u_1}-\sum_{j=2}^4\frac{(u_i-u_j)(u_k-u_j)}{2(u_1-u_j)} \frac{\partial^2 g(u^0)}{\partial u_1^2}\right) 
\eea after using properties (6.4) and (6.2). (Here $ \{i,j,k\}=\{2,3,4\}$). On the other hand, writing
\bea  C(u^0)=\frac{1}{2^4\delta^2(u^0)
}g(u^0)\prod_{i=1}^4\frac{\partial^2 g(u^i)}{\partial u_0^2} \eea as a sum of simple fractions in the variable  $u_1$ and using Lemma 6.2, one obtains the following:
\bea  C(u^0) = \frac{1}{\delta^2(u^{01})}
\left(\frac{1}{2}\frac{\partial g(u^0)}{\partial
u_1}-\sum_{j=2}^4\frac{(u_i-u_j)(u_k-u_j)}{2(u_1-u_j)}
\frac{\partial^2 g(u^0)}{\partial u_1^2} \right)^2+\frac{3}{4}= \eea
\bea = E_1(u^0)^2 +\frac{3}{4}.\eea
 Since this works equally well for each $E_i(u^0)$, we have proven formula (6.13) which in turn implies 
\bea \frac{G(u)}{\delta^2(u)} = \left( \sum_{i=1}^4\frac{g(u^0)}{\delta (u^{0i})(u_0-u_i)} + E(u^0) \right)^2 + \frac{3}{4}.\eea
 Thus formula (6.8) is proven. 

 We now know that for $w$ generic, the fiber $(\cU/\hat{G})_w$ has two irreducible components. Going back to $(\cS/\hat{G})_w$, notice that the two corresponding irreducible components must be disjoint. Indeed, formula (6.8) shows that the points of intersection of the above two components, if existent, should be among the van Geemen lines. This would contradict Lemma 3.1 which states that for $t\not=\mu\sqrt[5]{ \frac{2^7}{3}}$, the van Geemen lines define smooth points of $\cS_t$. Also, a simple check shows that the group $\hat{G}$ acts freely on the set of these lines and they also give smooth points of $(\cS/\hat{G})_w$ (see also Lemma 6.3). Thus for $w\not= 0, \frac{2^7}{3}$,  $(\cS/\hat{G})_w$ has (at least) two connected components. 
\end{proof}

 More information on the fiber $(\cS/\hat{G})_w$ can be derived by studying the action of the group $\hat{G}$
 on $\cS$. One can decompose this action into the action of $G$,
 which fixes the fiber, and that of $\mu=\prod_j\mu_j$ on
$t\in\CC$. The group of fifth roots of unity acts freely
on $(\CC)^*$. The smoothness of the generic $(\cS/G)_t$ follows via
 the following lemma:
 
\begin{lemma}
 Let the group G act on the projective space $\PP^4$ with coordinates
 $(x_0:...:x_4)$. Then the only lines in $\PP^4$ with nontrivial
 isotropy group are those inside the hyperplanes $H_i=(x_i=0)$. \end{lemma}

\begin{proposition}
  For t generic,  $\cH_t/G$ consists of two isomorphic smooth hyperelliptic genus 6 curves together with 15 isolated points.
\end{proposition}
\begin{proof}
  By the preceding lemma, $G$ acts freely on $\cS_t$ for all $t\not=0$. Thus the generic quotient $\cS_t/G$ consists of two smooth curves of genus 6, as may be computed by the Hurwitz formula. Each of the 375 isolated points has an isotropy group of order 5, thus one obtains 15 isolated points in the quotient.

 One can see that each component $C^1_t$, $C^2_t$ of the generic fiber
$\cS_t/G$ is hyperelliptic by examining the quotients of the Van
Geemen points. There are 5000 such points, organized in 10 families of
500. Each such family is the intersection of $\cS_t$ with a
hyperplane $(p_{ij}=0)$ in $\PP^9$ with Pl\"{u}cker coordinates. On
the other hand, formula (6.8) shows that the action of $\tau\in\ZZ_2$
on $\cS_t/G$, induced by the interchanging of two coordinates in
$\PP^4$, is mapping isomorphically one component into the other, so
the 500 points in each family split into two families of 250. Each of
these families corresponds to the linear divisor given by the
hyperplane section of the corresponding component in $\PP^9$. After
quotienting by the G-action, one obtains 10 families of 2 points, and
linear equivalence is preserved. The linear system $ |p_1+p_2|$
in $C^i_t$, where $p_1$, $p_2$ are the classes of the lines through:
\bea \left( \begin{array}{c} (1:\xi:\xi^2:0:0),\\
                             (1:1:1:a:b)\end{array}\right)
\mbox{ and } \left( \begin{array}{c} (1:\xi^2:\xi:0:0),\\
                             (1:1:1:b:a)\end{array}\right). \eea
 is base-point-free and gives a 2-1 morphism $q^i_t$ into
$\PP^1$. Moreover, since $a$ and $b$ above depend holomorphically on
$t$, it follows that the morphism $q^i_t$ does also. Notice also that
even if $C^i_t$ is singular, one may still repeat the above argument
for its normalization, since the Van Geemen points are smooth for
$t\not=\ts$. This even works for $t=\ts$, since the Van Geemen points
there are smooth points of the reduced fiber. Thus one establishes a
rational map from $\cS/G$ into $\PP^1$, 2-1 on fibers.
\end{proof}
 
\begin{remark}
 By Hurwitz formula, the irreducible components of $\cH_{0}/\hat{G}$ are all rational. Indeed, with the notations of section 1 and Theorem 5.2, there is an order 5 subgroup of $\hat{G}$ permuting the cones $(C_{i,j,\mu})_{\mu}$, another order 5 subgroup fixing each element of  $C_{i,j,\mu}$, and from the order 25 remaining action:
\bea 2\cdot(6-1)= 25\cdot 2(0-1) + 15\cdot 4.\eea
  Since there are 10 such irreducible components of $\cH_{0}/\hat{G}$, intersecting two by two at 15 points, one immediately gets a total genus of 6 for the reduced $\cH_{0}/\hat{G}$.
\end{remark}

\section{ Properties of the fibers $\cS_t$ }

\begin{proposition}
The fiber of the family $(\cS_t)_t$ at $t\not= 0, \ts$ consists of two isomorphic connected curves of genus 626 and degree 250 in G(2,5). The fibers at $t=0$ and $\ts$ are connected and are contained in the family with multiplicity 2. The set of all the connected curves in the family can be parametrized by a quasi-projective curve which extends to a genus 2 hyperelliptic curve $C_2$.
\end{proposition}
 
\begin{proof}

 Proposition 6.1 shows that for $w\not=0, \frac{2^7}{3}$, $(\cS/\hat{G})_w$ consists of two disjoint isomorphic connected curves, so the Stein factorization of the morphism 
\bea \begin{CD} (\cS/\hat{G}) @> {\pi/\hat{G}} >> (\CC/\hat{G}) \end{CD} \eea is\bea  \begin{CD} (\cS/\hat{G}) @> {\pi_2/\hat{G}} >> (C_2/\hat{G}) @> {\beta/\hat{G}} >> (\CC/\hat{G}) \end{CD}, \eea
 where $\beta/\hat{G}$ extends to a double cover of $\PP^1$ branched at 0 and $\frac{2^7}{3}$.

 A similar statement is true for the fibers $\cS_t$. The fact that the $\hat{G}$-action does not influence the number of connected components of the generic fiber may be checked in different ways: either by recalling the behavior of $\cH$ in a neighborhood of $\cH_0$ described in Section 2, or computationally be checking that the equations (6.6) are in general irreducible.

 Consider the closure $\bar{\cS}$ of the quasi-projective surface $\cS$ in $G(2,5)\times \PP^1$. Let  $\tilde{\cS}$ be the normalization of $\bar{\cS}$. By Stein factorization theorem, the morphism
 \bea \begin{CD}
 \tilde{\cS} @> \tilde{\pi} >> \PP^1 \end{CD} \eea
factors as  \bea \begin{CD}
 \tilde{\cS} @> {\tilde{\pi}_2}>> C_2 @>{\bar{\beta}} >> \PP^1, \end{CD} \eea
where $\tilde{\pi}_2$ is a projective morphism with connected fibers and $\bar{\beta}$ is a 2-1 morphism, ramified over $0$ and $\ts$, for $\mu^5=1$, while the curve $C_2$ must be a hyperelliptic curve of genus 2.
\end{proof}

\bigskip

\section{The stable limits $\tilde{\cS}_{\infty }$ and $\tilde{\cS}_{\infty }/\hat{G}$}

 Consider the following two compactifications of the surface $\cS$ at $t=\infty$:

 First, let  $\bar{\cS}$ be the closure of the quasi-projective surface $\cS \subset G(2,5)\times \CC$ in $G(2,5)\times \PP^1$. This surface is in fact not normal and the fiber $\bar{\cS}_{\infty}$ has some embedded points, as Lemma 8.2 will show.  

 Second, let $\tilde{G}(2,5)$ denote the blow-up of the Grassmannian $G(2,5)$ of lines in $\PP^4$, along the subvarieties $G_i(2,4)$ parametrizing lines inside the coordinate hyperplanes of $\PP^4$. Let $\tilde{\cS}$ be the closure of $\cS$ in $\tilde{G}(2,5)\times \PP^1$. The surface $\tilde{\cS}$ turns out to be smooth, as will be seen in Theorem 8.3. The fiber $\tilde{\cS}_{\infty}$ is the stable limit of the family $\cS$ at $t=\infty$ and consists of two connected components $\tilde{\cC}^{\xi}$ and $\tilde{\cC}^{\xi^2}$, which are isomorphic and reducible:
\bea \tilde{\cC}^{\xi}= \bigcup_{i=0}^4 \tilde{\cC}^{\xi}_i, \eea
with $\tilde{\cC}^{\xi}_i$ smooth isomorphic curves intersecting pairwise transversely at 25 points.

 The core of the argument here will consist in understanding the structure of the morphism $p: \tilde{G}(2,5) \to G(2,5)$ and the way the image of $\bar{\cS}$ in $G(2,5)$ transforms under $p$:

 Set $\{h,i,j,k,l\}=\{0,\ldots{},4\}.$ Working on the affine subsets $\AA^6\cong U_{ij} \subset G(2,5)$,  $U_{ij}=(p_{ij}\not= 0)$, with coordinates $(x_{k}, x_{l}, x_{h}, y_{k}, y_{l}, y_{h})$, one sees that blowing-up the ideal 
\bea (x_{k},y_{k})\cdot (x_{l},y_{l})\cdot (x_{h},y_{h}) \eea
of $\left( \bigcup_{i=0}^{4} G_i(2,4)\right) \cap  V_{ij} $ in $V_{ij}$ amounts to taking product of three separate blow-ups, in the directions $(x_{k},y_{k})$, $(x_{l},y_{l})$, and $(x_{h},y_{h})$:
\bea \begin{CD} (Bl_{(0,0)}\AA^2)\times (Bl_{(0,0)}\AA^2)\times (Bl_{(0,0)}\AA^2) @>>> \AA^6. \end{CD} \eea
 The fibers contained in the exceptional divisors are:

\begin{itemize}

\item $\PP^1$ over points of $G_{i}(2,4) \backslash
(\bigcup_{j\not= i} G_{ij}(2,3))$, where each $ G_{ij}(2,3)$ 
denotes the Grassmannian of lines in the
plane $(x_i=x_j=0)$

\item $\PP^1\times \PP^1$ over points of  $ G_{ij}(2,3)$ except
for those corresponding to the lines  $(x_i=x_j=x_k=0)$

\item $\PP^1\times \PP^1\times \PP^1$ over the points of $G(2,5)$
corresponding to lines  $(x_i=x_j=x_k=0)$.
\end{itemize}

  Furthermore, if one identifies $\bar{\cS}$, respectively $\tilde{\cS}$ with their images in $G(2,5)$ and $\tilde{G}(2,5)$, one will obtain that
$$ p^*(\bar{\cS}) = \tilde{\cS}  \cup \left( \bigcup_{i=0}^4 E_i \right)  $$
scheme-theoretically, where $E_i$ are the exceptional divisors in $\tilde{G}(2,5)$. This allows for a simple algebraic presentation of $\tilde{\cS}$, to the effect that $\tilde{\cS}$ is a smooth surface. 

  The action of the group $\hat{G}$ on $G(2,5)$ extends canonically to an action on  $\tilde{G}(2,5)$, making the morphism $$p: \tilde{G}(2,5) \to G(2,5)$$
$\hat{G}$-equivariant. A brief study of this action will show that
$\tilde{\cS}/\hat{G}$ is a smooth surface and the stable limit at $ w
=\infty$ of the family $\cS/\hat{G}$ is 
$$\tilde{\cS}_{\infty}/\hat{G} = (\tilde{\cC}^{\xi}/\hat{G}) \bigcup (\tilde{\cC}^{\xi^2}/\hat{G})$$ where $\tilde{\cC}^{\xi}/\hat{G}$ is the union of 5 smooth rational curves intersecting pairwise transversely.

\smallskip

 Here we will keep the notations introduced at the beginning of section 1. 

 In \cite{ak2}, the subset $I$ of G(2,5) consisting of all the lines
 incident to the 5 components of the base locus $\cB$ is described as
 a complete intersection surface in G(2,5) given by the following
 equations:\[ \sum _{ j\not=i; j=0}^4 p_{ij}^5=0 \] for i=0,4. One
 of these five equations is a linear combination of the others, which amounts to saying that whenever a line $l$ intersects 4 of the components $B_i$, it automatically intersects the fifth.

According to the possible multiplicities of intersection of a line $l$
with the various $B_i$-s, one distinguishes the following irreducible
components of $I$: 

 (1) 50 components, each of which consists of lines in a cone over one
     of the  $B_l$-s and having as vertex one of the points in
     $B_{ijk}$; 

 (2) 15 components, each of which consists of lines intersecting two
     of the curves  $B_{ij}$ and  $B_{kl}$; 

 (3) 10 components, each of which consists of lines intersecting one
     $B_{ij}$ and  $ B_k$ and $B_l$ and $B_h$ (but in general not
     $B_{kl}$, $B_{kh}$ or $B_{lh}$, nor $B_{ijk}$ $B_{ijl}$ or
     $B_{ijh}$);

 (4) The image of $\bar{\cS} \in G(2,5)\times\PP^1$ through projection
     on the Grassmannian $G(2,5)$.

 We will start with a description of $\bar{\cS}$ initiated by \cite{ak2}:

\begin{proposition}  $\bar{\cS}$ is isomorphic onto its image in $G(2,5)$, which is the closure of the set 
\bea \{ \{l\}\in G(2,5) / l\cap B_i\not=\emptyset,
l\cap B_{ij}=\emptyset, \forall i,j\in\{0,...4\}, i\not=j \}. \eea 
\end{proposition} 

\begin{proof}
 Indeed, if one considers the line $l$ embedded as $ \phi: \PP^1 \to \PP^4 $, then $\{l\}\in\cH$ if and only if the pull-back $ \phi^*f$ of the rational map on $\PP^4$ 
\bea f((x_0:...:x_4))=\frac{x_0^5+x_1^5+x_2^5+x_3^5+x_4^5}{x_0x_1x_2x_3x_4} 
\eea
is a constant map, i.e. $\phi^*(\cO_{\PP^4}(X_0-\sum_{i=0}^4 H_i)\cong
\cO_{\PP^1}$. Hence 
\bea  \{ \{l\}\in G(2,5) / l\cap B_i\not=\emptyset, 
l\cap B_{ij}=\emptyset, \forall i,j\in\{0,...4\}, i\not=j \} \subset \bar{\cS},\eea
 and on the other hand $\bar{\cS}$ is contained in the variety $I$ of
all the lines intersecting  the $B_i$-s. Since $\bar{\cS}$ is
irreducible, it is the closure of the above open set.
\end{proof}

\begin{lemma}
 The reduced structure of $\bar{\cS}_{\infty}$ is given by:
$$ \bar{\cS}_{\infty red} = \bigcup_{i=0}^4 (\cC^{\xi}_i \cup \cC^{\xi^2}_i)$$
where all $\cC^{\xi}_i$ and $\cC^{\xi^2}_i$ are smooth irreducible curves of genus 76, and they intersect pairwise transversely as follows: for each $i,j\in \{0,...,4\}$,  
$$\cC^{\xi}_i \cap \cC^{\xi^2}_i \cap \cC^{\xi}_j \cap \cC^{\xi^2}_j =\bigcup_{s=1}^{25}\{ P_{ij}^s\}.$$

 All the components $\cC^{\xi}_i$, $\cC^{\xi^2}_i$ are embedded with multiplicity 1 in $\bar{\cS}$ but $\bar{\cS}_{\infty}$ is not reduced, containing embedded points.

\end{lemma}

\begin{proof}
 Let  $G_i(2,4)$ denote the Grassmannian of lines contained in the
 hyperplane $H_i=(x_i=0)$ in $\PP^4$. By the description of $\bar{\cS}$ given above, the reduced structure of $\bar{\cS}_{\infty}$  consists of five copies of $ L_i $, the closure in $G_i(2,4)$ of the set
 \bea && \{ \{l\}\in G_i(2,4) /l\cap B_{ij}\not=\emptyset, 
 l\cap B_{ijk}=\emptyset, \forall j,k \in \{0,...4\}\backslash \{i\}, j\not=k \}. \eea  
 Take for example $i=4$. Again, $L_4$ is a subvariety of the complete intersection curve $C_4$ in G(2,4), which consists of all the lines in $(x_4=0)$ intersecting the curves $B_{4j}$.  $C_4$ is given by the following homogeneous ideal
\bea \cI=\left(\begin{array}{c} p_{01}^5+p_{02}^5+p_{03}^5, \\
  -p_{01}^5+p_{12}^5+p_{13}^5,\\
  -p_{02}^5-p_{12}^5+p_{23}^5,\\
  p_{01}p_{23}+p_{03}p_{12}+p_{02}p_{31}
\end{array}\right) \eea
The relation:
\bea  5p_{01}p_{02}p_{03}p_{12}p_{23}p_{31}( p_{03}p_{12}+\xi p_{02}p_{13})( p_{03}p_{12}+\xi^2 p_{02}p_{13}) = \eea
 \bea = p_{01}^5( -p_{02}^5-p_{12}^5+p_{23}^5)+p_{12}^5(p_{01}^5+p_{02}^5+p_{03}^5)-\eea
\bea - p_{02}^5( -p_{01}^5+p_{12}^5+p_{13}^5)
-( p_{01}^5p_{23}^5+(p_{03}p_{12}+p_{02}p_{31})^5)  \eea
 implies the existence of the following irreducible components of $C_4$, each appearing with multiplicity 1:

\begin{itemize}
\item 30 curves, corresponding to 30 cones in $\PP^3$, one through
each of the points of $B_{4ij}$, $ i,j \in \{0,...3\}$; (recall that
each $B_{4ij}$ is a union of 5 points $B_{4ij}^{\mu}$, one for each
fifth root of unity $\mu = $).  One example of the prime ideal of one of these curves would be:
\bea \cP_{034}^{\mu }=\left( p_{03}, p_{01}+\mu p_{02}, p_{31}+\mu p_{32},-p_{01}^5+p_{12}^5+p_{13}^5 \right) \eea
for $C_{034}^{\mu}$.

\item 2 other smooth irreducible curves:

$\cC_{4}^{\xi}$ given by the ideal:
\bea \left( \begin{array}{ccc} p_{03}p_{12}+\xi p_{02}p_{13}, & p_{01}p_{23}+\xi^2p_{02}p_{13}, & \\
   p_{01}^5+p_{02}^5+p_{03}^5,& 
  -p_{01}^5+p_{12}^5+p_{13}^5,&
  -p_{02}^5-p_{12}^5+p_{23}^5 \end{array} \right) \eea
 
$\cC_{4}^{\xi^2}$ given by the ideal:
\bea \left( \begin{array}{ccc} p_{03}p_{12}+\xi^2 p_{02}p_{13}, & p_{01}p_{23}+\xi p_{02}p_{13}, & \\
   p_{01}^5+p_{02}^5+p_{03}^5,& 
  -p_{01}^5+p_{12}^5+p_{13}^5,&
  -p_{02}^5-p_{12}^5+p_{23}^5 \end{array} \right) \eea
\end{itemize}

 These last two components will be part of $\bar{\cS}_{\infty red}$.

 By working with the affine cover of principal open sets \\
$U_{ij}=(p_{ij}\not= 0)$ of G(2,4), we can verify that these are all the irreducible components of $C_4$, that they are smooth, they appear with multiplicity 1 and the multiplicity of intersection of two components at a point of intersection is always 1. For example, on the affine set $U_{01}$, with coordinates 
$$\begin{array}{cc} x_i=\frac{-p_{1i}}{p_{01}},   & y_i=\frac{p_{0i}}{p_{01}} \end{array}$$
 for $i\in\{2,3\}$, the local equations of $C_4$ are:
\bea && 1+x_2^5+x_3^5=0\\
     && 1+y_2^5+y_3^5=0 \\
     && x_2^5-y_2^5+(x_2y_3-x_3y_2)^5=0 \eea
with
\bea  x_2^5-y_2^5+(x_2y_3-x_3y_2)^5= \eea
\bea = x_2^5(1+y_2^5+y_3^5)-y_2^5(1+x_2^5+x_3^5)+ \eea
\bea +5x_2x_3y_2y_3(x_2y_3+\xi x_3y_2)(x_2y_3+\xi^2 x_3y_2). \eea 

The components of $C_4$ intersect at the following points:
\begin{itemize}
\item $25\times 3=75$ points $l_{(ij),(hk)}^{\mu ,\nu }$, corresponding to lines passing through points $B_{ij4}^{\mu}$ and $B_{hk4}^{\nu}$, with $\{i,j,k,h\}=\{0,1,2,3\}$; each of these $l_{(ij),(hk)}^{\mu ,\nu }$ is at the intersection of exactly 2 components: those coming from cones in $\PP^3$ with vertices $B_{ij4}^{\mu}$, $B_{hk4}^{\nu}$, respectively, over the adequate quintic curves;

\item $25\times 4=100$ points $l_{ijk}^{\mu ,\nu }$, corresponding to lines inside the hyperplane $x_l=0$, passing through points $B_{ij4}^{\mu}$, $B_{ik4}^{\nu}$ and $B_{jk4}^{\mu / \nu}$ (here again $\{i,j,k,l\}=\{0,1,2,3\}$). Each of these 100 points lies at the intersection of 5 components of $C_4$: $C_{ij4}^{\mu}$, $C_{ik4}^{\nu}$, $C_{jk4}^{\mu / \nu}$, $\cC_{4}^{\xi}$ and $\cC_{4}^{\xi^2}$. 
\end{itemize}

 The above decomposition enables one to compute the Hilbert
polynomials $p_{\cC_{4}^{\xi}}(n)$ and $p_{S_{\infty}}(n)$ in
$G(2,5)\subset \PP^9$ from successive  sequences of the type
\bea 0 \to \cO_{\PP^9}/(\cI_1\cap \cI_2) \to \cO_{\PP^9}/ \cI_1 \oplus  \cO_{\PP^9}/ \cI_2 \to  \cO_{\PP^9}/(\cI_1 + \cI_2) \to 0:  \eea
\bea  p_{\cO_{\PP^9}/ \cI_1}(n)+p_{\cO_{\PP^9}/ \cI_2}(n) = p_{\cO_{\PP^9}/(\cI_1\bigcap \cI_2)}(n)+p_{\cO_{\PP^9}/(\cI_1 + \cI_2)}(n). \eea
Thus:
\bea \sum_{\{i,j\} \subset\{0,1,2,3\}, \mu}  p_{C_{ij4}^{\mu}}(n)+p_{\cC_{4}^{\xi}}(n)+p_{\cC_{4}^{\xi^2}}(n)=p_{C_4}(n) + 75 + \left(\begin{array}{c} 5\\
2\end{array}\right)100.\eea
\bea  p_{C_{ij4}^{\mu}}(n)=5n-5, \eea
since these curves embed in $\PP^9$ as plane quintics.
\bea p_{\cC_{4}^{\xi}}(n) = p_{\cC_{4}^{\xi^2}}(n). \eea 
\bea p_{C_4}(n) = 250n - 1375, \eea
since $C_4$ is a $(5,5,5,2)$ complete intersection in $\PP^5 \subset \PP^9$. Thus
\bea  p_{\cC_{4}^{\xi}}(n)=\frac{1}{2}\{250n-1375+75+1000-30(5n-5)\}=50n-75. \eea
\bea \bar{\cS}_{\infty red}= \bigcup_{i\in\{0,...4\}}\left(\cC_{i}^{\xi}\bigcup\cC_{i}^{\xi^2} \right) \eea
 and each point $l_{ijk}^{\mu ,\nu }$ is at the intersection of exactly 4 components ($\cC_{h}^{\xi}$,  $\cC_{h}^{\xi^2}$, $\cC_{l}^{\xi}$ and $\cC_{l}^{\xi^2}$, where $\{0,...,4\}\backslash\{i,j,k,l,h\}$). There are 250 such points when $\{i,j,k\}\subset\{0,...,4\}$.  These points will also be denoted by $\{P_{lh}^s\}_{l\not=h\in \{0,...,4\}}$, with $s \in\{1,\ldots{}25\}$. In our example $h=4$.

Notice also that none of the 375 points $l_{(ij),(hk)}^{\mu ,\nu }$ in the base locus is contained in $\bar{\cS}_{\infty red}$. Thus
\bea  p_{\bar{\cS}_{\infty red}}(n)=\sum_{i\in\{0,...4\}}2p_{\cC_{i}^{\xi}}(n)-\left(\begin{array}{c} 4\\
2\end{array}\right)250 =  500n-2250.  \eea

 Comparison with the formula for the Hilbert polynomial  of the family $(\cS_t)_t$  in Proposition 4.4 shows that the fiber $\bar{\cS}_{\infty}$ must be embedded with multiplicity 1 inside the family $(\cS_t)_t$. It also shows that there must be some embedded points in $\bar{\cS}_{\infty}$, which will give singular points of $\bar{\cS}$. We will see shortly that these points are exactly $\{P_{lh}^s\}$
, thus each appearing with multiplicity 5 in $\bar{\cS}_{\infty}$.

\end{proof}

 We now proceed with the description of the strict transform $\tilde{\cS}$ of $\bar{\cS}$:

\begin{theorem}
 The surface $\tilde{\cS}$ is smooth.

 The fiber $\tilde{\cS}_{\infty}$ is the stable limit of the family $\cS$ at $t=\infty$ and consists of two connected components $\tilde{\cC}^{\xi}$ and $\tilde{\cC}^{\xi^2}$, which are isomorphic and reducible:
\bea \tilde{\cC}^{\xi}= \bigcup_{i=0}^4 \tilde{\cC}^{\xi}_i \eea
with $\tilde{\cC}^{\xi}_i$ smooth isomorphic curves of genus 76, intersecting pairwise transversely at 25 points each.
\end{theorem}

\begin{proof}
 By examining the morphism $p: \tilde{G}(2,5) \to G(2,5)$ locally, we will see that no fibers of the exceptional divisors are contained in $\tilde{\cS}$, that for each pair $\cC^{\xi}_i $ and $\cC^{\xi^2}_i $, the two components are separated and $p$ desingularizes $\bar{\cS}$.
 
  We will work over $\AA^6\cong U_{01} \in G(2,5)$. Looking at the restriction of $p$:
\bea \begin{CD} (Bl_{(0,0)}\AA^2)^3 @>>> \AA^6 \end{CD} \eea
  one finds two types of affine charts, which after all possible permutations cover the entire $\tilde{V}_{01}$:\\
(1) $\AA^6 $ with coordinates $\{u_i,y_i\}_{i \in \{2,3,4\}}$, where the original coordinates \\
 $\{x_i,y_i\}_{i \in \{2,3,4\}}$ on $U_{01}$ satisfy:
$$ x_i=u_iy_i.$$
 The exceptional divisor $E_i=(y_i=0)$ and $u_i$ is the coordinate along the fiber of the exceptional divisor over the corresponding component of the blow-up locus.\\
(2)  $\AA^6 $ with coordinates $\{u_2,y_2, u_3, y_3, x_4, v_4\}$.
$$x_2=u_2y_2 \mbox{ , }  x_3=u_3y_3 \mbox{ , } y_4=x_4v_4. $$
 $(y_2)$, $ (y_3)$ and $(x_4)$ give the exceptional divisors.

 Case (1) is closely related to the setup of Proposition 4.1. Indeed, in this local picture formula (4.1) becomes 
\bea \phi(\alpha :\beta ) = (\alpha :\beta :y_2(u_2\alpha +\beta ):y_3(u_3\alpha +\beta ):y_4(u_4\alpha +\beta )) \eea
 and thus by the same reasoning as in Proposition 4.1:
\bea \left( \begin{array}{ccc}
                          u_2 & u_3 & u_4 \\
                          u_2^2 & u_3^2 & u_4^2 \\
                          u_2^3 & u_3^3 & u_4^3 \\
                          u_2^4 & u_3^4 & u_4^4 
 \end{array} \right)  \left( \begin{array}{c} y_2^5 \\ y_3^5 \\ y_4^5 \end{array} \right) =
 5t\prod_{j=2}^4 y_j \left( \begin{array}{c}  \frac{1}{5}s_0(u) \\ \frac{1}{10}s_1(u) \\ \frac{1}{10}s_2(u) \\ \frac{1}{5}s_3(u) \end{array} \right) \eea 
plus the separate equations 
 \bea  \sigma_5+1=\sigma_0+1=0, \eea
where the notations are the same as in section 4, except that in this case 
\bea u=(u_2:u_3:u_4) \eea
As before, one sees that the fourth row of the augmented matrix is just a linear combination of the others. After removing the exceptional divisors $E_2$, $E_3$, $E_4$, one is left with the irreducible surface $\tilde{\cS}$ given by:
\bean \sigma_5+1=\sigma_0+1= \sigma_2 - \frac{1}{2}s_1 \sigma_1= \sigma_3 - \frac{1}{2}s_2\sigma_1=0. \eean
  Consider the intersection with the exceptional divisor $E_4$: Then $y_4=0$ and $\sigma_1, \sigma_2, \sigma_3$ satisfy in this case the identity:
$$ \si_3 -(u_2+u_3)\si_2+u_2u_3\si_1 = 0 $$
which in terms of the equations (8.1) is equivalent to:
\bean \si_1 (u_2^2-u_2u_3+u_3^2) = 0  \eean
 The conditions $u_2= -\xi u_3$ and $u_2= -\xi^2 u_3$, where $\xi^2+\xi+1=0$,  give the direct transforms $\tilde{\cC}^{\xi}_4$ and $\tilde{\cC}^{\xi^2}_4$ of $\cC^{\xi}_4$ and $\cC^{\xi^2}_4$. Notice that in this open set  $\tilde{\cC}^{\xi}_4$ and $\tilde{\cC}^{\xi^2}_4$ do not intersect, because condition $\sigma_5+1=0$ excludes the case $u_2=u_3=0$. Also, in this open set $\tilde{\cS}$ contains no fibers in the exceptional divisor $E_4$: such fibers would satisfy $\si_1=0$, which would also imply $\si_2=0$. But these two equations in $u_2,u_3,y_2,y_3$ have no common solutions with $\sigma_5+1=\sigma_0+1= 0.$
 
 The intersection $\tilde{\cS} \cap E_4 \cap E_3$ immediately gives:  $$\begin{array}{ccc}y_3=y_4=0, &y_2=-\mu, &u_2=-\nu\end{array}$$  for two fifth roots of unity $\mu$ and $\nu$. Also, $\sigma_2 - \frac{1}{2}s_1 = 0 $ becomes $$u_2-u_3-u_4=0,$$ so in view of equations (8.1) and (8.2), one sees that $\tilde{\cC}^{\xi}_3$ intersects $\tilde{\cC}^{\xi}_4$ at 25 points given by $$\begin{array}{ccccc}y_3=y_4=0, & y_2=-\mu, & u_2=-\nu, & u_3=\xi^2\nu, &  u_4=\xi\nu \end{array}   $$
and  $\tilde{\cC}^{\xi^2}_3$ intersects $\tilde{\cC}^{\xi^2}_4$ at 25 points given by $$\begin{array}{ccccc}y_3=y_4=0, & y_2=-\mu, & u_2=-\nu, & u_3=\xi\nu, & u_4=\xi^2\nu \end{array} $$
 and these are all the points in  $\tilde{\cS} \cap E_4 \cap E_3$.

  In the above, $\tilde{\cC}^{\xi}_4$ and $\tilde{\cC}^{\xi^2}_4$ have
 been fixed; after permutations, the appropriate choices of notation for $\tilde{\cC}^{\xi}_i$ and $\tilde{\cC}^{\xi^2}_i$ can be made, by keeping track of the signatures of the permutations. In this way the set of components of $\tilde{\cS}_{\infty}$ found thus far splits in two subsets  $\{\tilde{\cC}^{\xi}_i\}_i$ and  $\{\tilde{\cC}^{\xi^2}_i\}_i$, such that all the curves in one subset intersect pairwise transversely, but do not intersect the components in the other subset.
 The smoothness of $\tilde{\cS}$ along the points of $\tilde{\cS}_{\infty}$ is easily checked by differentiating equations (8.2) and using the conditions derived from (8.3).  
 
 The only thing that remains to be checked is existence of fibers in the exceptional divisors, or intersection points for $\tilde{\cC}^{\xi}_i$ and $\tilde{\cC}^{\xi^2}_i$ in sets of type (2). Here equations (8.2) are written for
$$\si_i= u_2^iy_2^5 + u_3^iy_3^5 + x_4^5v_4^i, $$
$$s_0=v_4, s_1=(u_3+u_3)v_4+1, s_2=u_2u_3v_4+u_2+u_3, s_3=u_2u_3.$$
 A simple check following the lines of the previous computations shows that there are no extra intersection points, nor fibers in the exceptional divisors.
\end{proof}

 The genus computation for $\tilde{\cS}_{\infty}$ gives exactly:
$$1250= 2\cdot (5\cdot 75 +250) $$

\begin{theorem}
 The action of the group $\hat{G}$ on $G(2,5)$ extends canonically to an action on  $\tilde{G}(2,5)$, making the morphism 
$$p: \tilde{G}(2,5) \to G(2,5)$$
$\hat{G}$-equivariant. $\tilde{\cS}/\hat{G}$ is a smooth surface. The stable limit at $ w =\infty$ of the family $\cS/\hat{G}$ is 
$$\tilde{\cS}_{\infty}/\hat{G} = (\tilde{\cC}^{\xi}/\hat{G}) \bigcup (\tilde{\cC}^{\xi^2}/\hat{G})$$ where $\tilde{\cC}^{\xi}/\hat{G}$ is the union of 5 smooth rational curves intersecting pairwise transversely.
\end{theorem}

\begin{proof}
 Working locally over $\AA^6\cong U_{01} \in G(2,5)$, and considering for $\tilde{G}(2,5)$ the local coordiantes set up in the proof of Theorem 8.3, it is not difficult to derive the action of 
$\hat{G}$ on $\tilde{G}(2,5)$: If \bea g=(1:\mu_1:\mu_2:\mu_3:\mu_4) \in \hat{G} \eea 
and $\{l\} \in  G(2,5)$ is the line passing through $(1:0:x_2:x_3:x_4)$, $(0:1:y_2:y_3:y_4)$, then\bea g^*(x_{i})=\mu_ix_{i}, \eea
\bea g^*(y_{i})=\mu_iy_{i}/\mu_1, \eea
and since $x_{i}= u_{i}y_{i}$, it follows that
\bea g^*(u_{i})=\mu_1u_{i} \eea
for all $i\in\{2,3,4\}.$

 Recall from  Theorem 5.2 that $\cS/\hat{G}$ is smooth. The proof that all points of $\tilde{\cS}_{\infty}/\hat{G}$ are smooth in $\tilde{\cS}/\hat{G}$ follows the same argument as in the proof of Theorem 5.2. A general point of $\tilde{\cC}^{\xi}_4$, for example, has an order 5 isotropy group, generated by an element of the form $(1:1:1:1:\mu)$, with $\mu^{5}=1$. A point in  $\tilde{\cC}^{\xi}_4 \cap \tilde{\cC}^{\xi}_3$ has an order 25 isotropy group, generated by two elements of the form $(1:1:1:1:\mu)$, $(1:1:1:\mu:1)$, with $\mu^{5}=1$. Each of these elements of $g$ gives a pseudoreflection of $\tilde{\cS}$, because the set of points in $\tilde{\cS}$ fixed by it is the divisor  $\tilde{\cC}^{\xi}_4 \cup \tilde{\cC}^{\xi^{2}}_4$,  $\tilde{\cC}^{\xi}_3 \cup \tilde{\cC}^{\xi^{2}}_3$ respectively. Again following \cite{st},  $\tilde{\cS}/\hat{G}$ is smooth.

 The Hurwitz formula for the genus $g'$ of $\tilde{\cC}^{\xi}_{4}/\hat{G}$:
\bea   
2(76-1) = 125 \cdot (g'-1) + 100\cdot 4\eea
implies $g'=0$. The 25 points of  $\tilde{\cC}^{\xi}_4 \cap \tilde{\cC}^{\xi}_3$ are all in the same orbit, thus giving one point in the quotient.
Notice that the above numbers sum up to a total genus of 6 for  $\tilde{\cC}^{\xi}/\hat{G}$.
\end{proof}

\bigskip

\section{The class $[\bar{\cS}]$ on $G(2,5)$}

 Let $F(1,2,5)$ denote the flag variety
\bea F(1,2,5) = \{ (x,[l]) / [l]\in G(2,5), x\in l \} \eea with the natural morphisms:
\bea \begin{CD}  F(1,2,5) @ > p >> \PP^4 \\
 @V q VV \\
 G(2,5) \end{CD} \eea
and consider the induced diagram
\bea \begin{CD}  \bar{\cU} @ > {p|_{\cU}} >> V \\
 @V {q|_{\bar{\cU}}} VV \\
 \bar{\cS} \end{CD} \eea
where $\bar{\cU}$ is the pull-back of $F(1,2,5)$ on $\bar{\cS} $ and $V = p(\cU)$.

 Here we will use the incidence variety $I$ and the description of its components from Section 8 to compute the class $[\bar{\cS}]$ of $\bar{\cS}$ in the cohomology ring of the Grassmannian, and the degree of $V$ as threefold in $\PP^4$.
  The method is to compute the class of each component of $I$ and then eliminate the irrelevant components to get $\bar{\cS}$. Recall that $I$ comprises:\\
 (1) 50 components, each of which consists of lines in a cone over one of the  $B_l$ -s and having as vertex one of the points in $B_{ijk}$, \\
 (2) 15 components, each of which consists of lines intersecting two of the curves  $B_{ij}$ and  $B_{kl}$, \\
 (3) 10 components, each of which consists of lines intersecting one  $B_{ij}$ and  $ B_k$ and $B_l$ and $B_h$ (but in general not  $B_{kl}$, $B_{kh}$ or $B_{lh}$, nor $B_{ijk}$ $B_{ijl}$ or $B_{ijh}$),
 (4) $\bar{\cS}$.

\begin{lemma}
 Let $I$ be given in the Grassmannian $G(2,5)$ by the ideal 
$$\cI = ( \sum_{j=0}^4 p_{ij}^5 )_{i\in\{1,\ldots{},4\}}.$$
 Then all the 2-dimensional components of $I$  appear with multiplicity 1 in the ideal $\cI$.  
\end{lemma}

\begin{proof}

  The proof of Proposition 8.1 implies the present lemma for those components of $I$ contained in one of the Grassmannians $G_i(2,4)$. In particular, $\bar{\cS}$ appears with multiplicity 1. Consider now the ideal $\cI$ localized on $U_{01}=(p_{01}\not=0)$, with the usual coordinates $\{x_i,y_i\}_{i\in\{2,3,4\}}$. The generators of $\cI$:
$$ f_1=1+x_2^5+x_3^5+x_4^5, $$
$$ f_2=1+y_2^5+y_3^5+y_4^5, $$
$$ f_3=x_3^5-y_3^5+(x_3y_2-x_2y_3)^5+(x_3y_4-x_4y_3)^5,$$
$$f_4=x_4^5y_0^5-x_1^5y_4^5+(x_4y_2-x_2y_4)^5+(x_4y_3-x_3y_4)^5.$$

satisfy the following relations:
\bean &&  f_3=5x_3y_3g_3+x_3^5f_2-y_3^5f_1 \\
      &&  f_4=5x_4y_4g_4+x_4^5f_2-y_4^5f_1 \eean
where  for $\{i,j,k\}=\{2,3,4\}$, $$g_{j}= p_{ij} p_{ij}^{\xi} p_{ij}^{\xi ^{2}}x_iy_i-p_{jk} p_{jk}^{\xi} p_{jk}^{\xi ^{2}}x_ky_k,$$
$$ \begin{array}{ccc} p_{ij}=x_iy_j-x_jy_i, &
 p_{ij}^{\xi}=x_iy_j+\xi x_jy_i, &
 p_{ij}^{\xi ^{-1}}=x_iy_j+\xi ^{-1}x_jy_i \end{array},$$
with $\xi^2+\xi +1=0$.
 Relations 9.1 and 9.2 highlight the existence of 2-dimensional components of the type:
$$\begin{array}{ccc} (x_i, y_j, f_1, f_2), & (x_i, x_j, f_1, f_2), & (y_i, y_j, f_1, f_2), \end{array}$$
$$ \begin{array}{cc} (p_{ij}, p_{jk}, p_{ki}, f_1, f_2), & (p_{ij}, p_{jk}^{\xi}, p_{ki}^{\xi ^{2}}, f_1, f_2) \end{array} $$ 
for $\{i,j,k\}=\{2,3,4\}$
  A straighforward exercise shows that these components also appear with multiplicity one and they cover all types (1)-(3) of components written above.
\end{proof}

For computations in the cohomology ring of $G(2,5)$, we will use the
following Schubert cycles: $\si_{(a_1,a_2)}$ is the class of $\{l | \dim (l\cap V_{3-a_i+i}\geq i \}$, for a fixed flag $V_1\subset V_2 \subset \PP^4 $.

\begin{proposition}
 $$ [\bar{\cS}] = 5^3\cdot 3\si_{22} +  5^3\cdot 2\si_{31}. $$
\end{proposition}

\begin{proof}
The different types of components of $I$ mentioned above give the
following cycles:\\
(0) $[I]$ corresponds to the cycle of lines intersecting 4 quintic surfaces in $\PP^4$:
$$ [I] =5^4 \si_1^4 =5^4\cdot 2 \si_{22} +5^4\cdot 3; \si_{31} $$
(1) $[I_1]$ corresponds to the cycle of lines contained in a cone over a quintic surface:
$$ [I_1]= 5 \si_{31}; $$
(2) $[I_2]$ corresponds to the cycle of lines intersecting two plane quintic curves:
$$ [I_2]= 25 \si_{20}^2 = 25( \si_{22} +\si_{31});  $$
(3) the cycle of lines intersecting a degree 5 plane  curve and two degree 5 surfaces consists of $[I_3]$ together with 3 copies of $[I_2]$ and 15 copies $[I_1]$ (again, coming from the geometry of the base locus $B$). Thus:
\bea [I_3] & = & 125 \si_{10}^2\cdot \si_{20} - 3\cdot 25( \si_{22} +\si_{31}) - 15\cdot 5 \si_{31}\\
  & = & 125 (\si_{11} +\si_{20})\si_{20} - 75( \si_{22} +\si_{31})-75 \si_{31} =\\  & = & 125 (2 \si_{31}+\si_{22})- 75( \si_{22} +\si_{31})-75 \si_{31} =\\
  & = & 100  \si_{31} + 50 \si_{22} \eea
Finally, 
$$ [\bar{\cS}] = [I] - 50 [I_1] -15 [I_2] -10 [I_3] =  5^3\cdot 3\si_{22} +  5^3\cdot 2\si_{31}. $$
 The degree of the surface $\bar{\cS}$ in $\PP^9$ is thus:
$$ deg (\bar{\cS}) = \int_{G(2,5)} \si_1^2 \cdot [\bar{\cS}] = 5^4. $$
\end{proof}

\begin{proposition}
The locus V covered by all the lines inside the Dwork pencil is a degree 250 threefold in $\PP^4$. 
\end{proposition}

\begin{proof}
 Recall the morphisms
\bea \begin{CD}  F(1,2,5) @ > p >> \PP^4 \\
 @V q VV \\
 G(2,5) \end{CD} \eea
 Let $H = p^*\cO_{\PP^4}(1)$. The degree of $V$ is:
$$deg(V) = \int_{F(1,2,5)} q^* [\bar{\cS}]\cdot H^3. $$
  By projection formula and the identity $q_*(H^3)=\si_2$:
 $$deg(V) = ( 5^3\cdot 3\si_{22} +  5^3\cdot 2\si_{31})\cdot \si_2 = 250.$$
\end{proof}

\begin{proposition}

For each $i \in \{0,...,4\}$, there is a finite, 30-1 morphism 
\bea\begin{CD}\tilde{\cS} @>>> B_i \end{CD}\eea
having 1-dimensional fibers over points of the type 
$$\begin{array}{cc} (1:-1:0:0:0), &(1:\xi:\xi^2:0:0) \end{array}$$ and their orbits through the actions of $\hat{G}$ and of the group of symmetries $S_5$, and branched over points of type $$(x_0:x_1:x_2:0:0)$$ with $x_0^5+x_1^5+x_2^5=0$ (and their orbits through the actions of  $\hat{G}$ and  $S_5$).

\end{proposition}

\begin{proof}

There is a well defined morphism from the blow-up of $G(2,5)$ along
 $G_i(2,4)$ to the hyperplane $H_i$ in $\PP^4$:if we think of this
 blow-up as sitting inside $F(1,2,5)$, then the morphism sends the pair $(x,[l]) \in F(1,2,5)$ to $l\cap H_i$ if $l$ is not contained in $H_i$, or to $x$ if $l \subset H_i$. The image of the strict transform of $\bar{\cS}$ through this morphism is exactly the Fermat surface $B_i$.
 Hence the morphism
\bea\begin{CD}\tilde{\cS} @>{\varphi_i}>> B_i. \end{CD}\eea
 To compute the degree of $\varphi_i$, set $i=0$ and consider a
 generic point $x$ of $B_0$. Consider the projection of vertex $x$
 down to $H_4$ and count the points where the images of the surfaces $B_1$ and
 $B_2$ intersect $B_4$. There are 125 such points: the projections of
 the 125 lines having contact with each of $B_0$,  $B_1$, $B_2$ and
 $B_4$ and therefore with $B_3$ also. Now we count how many of these
 lines are in components of type (1)-(3) of $I$: the 6 curves $B_{ij}$
 with $i,j \not=0$ give $150$ lines through $x$, that intersect
 $B_{ij}$ and the other $B_k$-s.  However, the lines through the 20 points of $B_{ijk}$, ($i,j,k \not=0$), have been counted three times in this process. By adjusting:
  $125-150+2\cdot20 =15$ lines.

 There is another special set of lines, which have been counted twice: lines through $x$ and intersecting two curves $B_{ij}$ and $B_{kl}$ with $\{i,j,k,l\}=\{1,2,3,4\}$. Indeed, lines intersecting $B_{ij}$ and $B_{kl}$ automatically intersect $B_0$ and by dimension count, the points of intersection with $B_0$ of all such lines cover the whole $B_0$. For each choice of pairs $(i,j)$ and $(k,l)$, there are 5 such lines passing trough $x$, since the projections of the planes $(x_i=x_j=0)$ and $(x_k=x_l=0)$ to the hyperplane $H_4$ intersect each other and $B_4$ in 5 points. A total of $15$ lines has thus been counted twice. Adding up, one obtains 30 lines through $x$ which are in the universal family $\cU$.

 The special fibers of the morphism $\varphi_i$ occur at the special
 lines studied: those in $\cS_0$ and the Van Geemen lines. 
\end{proof}

 Let $V_t=p_*q^*(\bar{\cS}_t)$, for any $t\in \PP^1$. $V_t$ is a degree 500 surface in $\PP^4$, as can be easily checked either at $t=0$ or $t=\infty$. 

 Notice that by our computations, $V \cap X_t = V_t \cup (\bigcup_{i=0}^4 B_i) $ is a reduced surface of degree $250\cdot 5 =1250$. As expected, each of the 5 surfaces $B_i$ appears with multiplicity 30 in $V \cap X_t$:
$$ 1250 = 500 + 30 \cdot 5 \cdot 5. $$
 
\begin{remark} In unpublished work, G.Pacienza has studied the\\
 2-dimensional locus of points where lines intersect a generic quintic threefold $X$ with multiplicity at least 5. He has shown that this locus has degree 650. If $X=X_{t
}$ is one of the quintics in the Dwork pencil, one can now state that the above locus is reducible to the degree 500 surface  $V_t$ and the 5 components of the base locus of $\cX$, each appearing with multiplicity 6. This is also in agreement with the numerical results of section 9.
\end{remark}


\providecommand{\bysame}{\leavevmode\hbox to3em{\hrulefill}\thinspace}

\end{document}